# LazySVD: Even Faster SVD Decomposition Yet Without Agonizing Pain

Zeyuan Allen-Zhu
zeyuan@csail.mit.edu
Institute for Advanced Study
& Princeton University

Yuanzhi Li
yuanzhil@cs.princeton.edu
Princeton University

**Abstract**

We study $k$-SVD that is to obtain the first $k$ singular vectors of a matrix $A$. Recently, a few breakthroughs have been discovered on $k$-SVD: Musco and Musco [19] proved the first gap-free convergence result using the block Krylov method, Shamir [21] discovered the first variance-reduction stochastic method, and Bhojanapalli et al. [7] provided the fastest $O(\mathsf{nnz}(A) + \mathsf{poly}(1/\varepsilon))$-time algorithm using alternating minimization.

In this paper, we put forward a new and simple LazySVD framework to improve the above breakthroughs. This framework leads to a faster gap-free method outperforming [19], and the first accelerated *and* stochastic method outperforming [21]. In the $O(\mathsf{nnz}(A) + \mathsf{poly}(1/\varepsilon))$ running-time regime, LazySVD outperforms [7] in certain parameter regimes without even using alternating minimization.

## 1 Introduction

The singular value decomposition (SVD) of a rank-$r$ matrix $A \in \mathbb{R}^{d\times n}$ corresponds to decomposing $A = V\Sigma U^\top$ where $V \in \mathbb{R}^{d\times r}$, $U \in \mathbb{R}^{n\times r}$ are two column orthonormal matrices, and $\Sigma = \mathrm{diag}\{\sigma_1, \dots, \sigma_r\} \in \mathbb{R}^{r\times r}$ is a non-negative diagonal matrix with $\sigma_1 \geq \sigma_2 \geq \cdots \geq \sigma_r \geq 0$. The columns of $V$ (resp. $U$) are called the left (resp. right) singular vectors of $A$ and the diagonal entries of $\Sigma$ are called the singular values of $A$. SVD is one of the most fundamental tools used in machine learning, computer vision, statistics, and operations research, and is essentially equivalent to principal component analysis (PCA) up to column averaging.

A rank $k$ partial SVD, or $k$-SVD for short, is to find the top $k$ left singular vectors of $A$, or equivalently, the first $k$ columns of $V$. Denoting by $V_k \in \mathbb{R}^{d\times k}$ the first $k$ columns of $V$, and $U_k$ the first $k$ columns of $U$, one can define $A_k^* \stackrel{\text{def}}{=} V_k V_k^\top A = V_k \Sigma_k U_k^\top$ where $\Sigma_k = \mathrm{diag}\{\sigma_1, \dots, \sigma_k\}$. Under this notation, $A_k^*$ is the the best rank-$k$ approximation of matrix $A$ in terms of minimizing $\|A - A_k\|$ among all rank $k$ matrices $A_k$. Here, the norm can be any Schatten-$q$ norm for $q \in [1, \infty]$, including spectral norm ($q = \infty$) and Frobenius norm ($q = 2$), therefore making $k$-SVD a very powerful tool for information retrieval, data de-noising, or even data compression.

Traditional algorithms to compute SVD essentially run in time $O(nd \min\{d, n\})$, which is usually very expensive for big-data scenarios. As for $k$-SVD, defining $\mathsf{gap} \stackrel{\text{def}}{=} (\sigma_k - \sigma_{k+1})/(\sigma_k)$ to be the relative $k$-th eigengap of matrix $A$, the famous subspace power method or block Krylov method [14] solves $k$-SVD in time $O(\mathsf{gap}^{-1} \cdot k \cdot \mathsf{nnz}(A) \cdot \log(1/\varepsilon))$ or $O(\mathsf{gap}^{-0.5} \cdot k \cdot \mathsf{nnz}(A) \cdot \log(1/\varepsilon))$ respectively if ignoring lower-order terms. Here, $\mathsf{nnz}(A)$ is the number of non-zero elements in matrix $A$, and the more precise running times are stated in Table 1.

| Paper | Running time (× for being outperformed) | GF? | Stoc? | Acc? |
|---|---|---|---|---|
| subspace PM [19] | $\widetilde{O}\big(\frac{k\mathsf{nnz}(A)}{\varepsilon}+\frac{k^2d}{\varepsilon}\big)$ × | yes | no | no |
| | $\widetilde{O}\big(\frac{k\mathsf{nnz}(A)}{\mathsf{gap}}+\frac{k^2d}{\mathsf{gap}}\big)$ × | no | | |
| block Krylov [19] | $\widetilde{O}\big(\frac{k\mathsf{nnz}(A)}{\varepsilon^{1/2}}+\frac{k^2d}{\varepsilon}+\frac{k^3}{\varepsilon^{3/2}}\big)$ × | yes | no | yes |
| | $\widetilde{O}\big(\frac{k\mathsf{nnz}(A)}{\mathsf{gap}^{1/2}}+\frac{k^2d}{\mathsf{gap}}+\frac{k^3}{\mathsf{gap}^{3/2}}\big)$ × | no | | |
| `LazySVD` **Corollary 4.3 and 4.4** | $\widetilde{O}\big(\frac{k\mathsf{nnz}(A)}{\varepsilon^{1/2}}+\frac{k^2d}{\varepsilon^{1/2}}\big)$ | yes | no | yes |
| | $\widetilde{O}\big(\frac{k\mathsf{nnz}(A)}{\mathsf{gap}^{1/2}}+\frac{k^2d}{\mathsf{gap}^{1/2}}\big)$ | no | | |
| Shamir [21] | $\widetilde{O}\big(knd+\frac{k^4d}{\sigma_k^4\mathsf{gap}^2}\big)$ (local convergence only) × | no | yes | no |
| `LazySVD` **Corollary 4.3 and 4.4** | $\widetilde{O}\big(knd+\frac{kn^{3/4}d}{\sigma_k^{1/2}\varepsilon^{1/2}}\big)$ $\Big($always $\leq \widetilde{O}\big(knd+\frac{kd}{\sigma_k^2\varepsilon^2}\big)\Big)$ | yes | yes | yes |
| | $\widetilde{O}\big(knd+\frac{kn^{3/4}d}{\sigma_k^{1/2}\mathsf{gap}^{1/2}}\big)$ $\Big($always $\leq \widetilde{O}\big(knd+\frac{kd}{\sigma_k^2\mathsf{gap}^2}\big)\Big)$ | no | | |
| All GF results above provide $(1+\varepsilon)\|\Delta\|_2$ spectral and $(1+\varepsilon)\|\Delta\|_F$ Frobenius guarantees | | | | |

Table 1: Performance comparison among direct methods. Define $\mathsf{gap} = (\sigma_k - \sigma_{k+1})/\sigma_k \in [0, 1]$. GF = Gap Free; Stoc = Stochastic; Acc = Accelerted. Stochastic results in this table are assuming $\|a_i\|_2 \leq 1$ following (1.1).

Recently, there are breakthroughs to compute $k$-SVD faster, from three distinct perspectives.

The *first breakthrough* is the work of Musco and Musco [19] for proving a running time for $k$-SVD that does not depend on singular value gaps (or any other properties) of $A$. As highlighted in [19], providing ***gap-free*** results was an open question for decades and is a more reliable goal for practical purposes. Specifically, they proved that the block Krylov method converges in time $\widetilde{O}\big(\frac{k\mathsf{nnz}(A)}{\varepsilon^{1/2}}+\frac{k^2d}{\varepsilon}+\frac{k^3}{\varepsilon^{3/2}}\big)$, where $\varepsilon$ is the multiplicative approximation error.[1]

The *second breakthrough* is the work of Shamir [21] for providing a fast ***stochastic*** $k$-SVD algorithm. In a stochastic setting, one assumes[2]

$$A \text{ is given in form } AA^\top = \frac{1}{n}\sum_{i=1}^{n} a_i a_i^\top \text{ and each } a_i \in \mathbb{R}^d \text{ has Euclidean norm at most } 1 \ . \quad (1.1)$$

Instead of repeatedly multiplying matrix $AA^\top$ to a vector in the (subspace) power method, Shamir proposed to use a random rank-1 copy $a_i a_i^\top$ to approximate such multiplications. When equipped with very ad-hoc variance-reduction techniques, Shamir showed that the algorithm has a better (local) performance than power method (see Table 1). Unfortunately, Shamir's result is (1) not gap-free; (2) not accelerated (i.e., does not match the $\mathsf{gap}^{-0.5}$ dependence comparing to block Krylov); and (3) requires a very accurate warm-start that in principle can take a very long time to compute.

The *third breakthrough* is in obtaining running times of the form $\widetilde{O}(\mathsf{nnz}(A) + \mathsf{poly}(k, 1/\varepsilon)\cdot(n + d))$ [7, 8], see Table 2. We call them ***NNZ results***. To obtain NNZ results, one needs sub-sampling on the matrix and this incurs a poor dependence on $\varepsilon$. For this reason, the polynomial dependence on $1/\varepsilon$ is usually considered as the most important factor. In 2015, Bhojanapalli et al. [7] obtained a $1/\varepsilon^2$-rate NNZ result using alternating minimization. Since $1/\varepsilon^2$ also shows up in the sampling complexity, we believe the quadratic dependence on $\varepsilon$ is tight among NNZ types of algorithms.

[1]In this paper, we use $\widetilde{O}$ notations to hide possible logarithmic factors on $1/\mathsf{gap}, 1/\varepsilon, n, d, k$ and potentially also on $\sigma_1/\sigma_{k+1}$.

[2]This normalization follows the tradition of stochastic $k$-SVD or 1-SVD literatures [12, 20, 21] in order to state results more cleanly.

All the cited results rely on ad-hoc non-convex optimization techniques together with matrix algebra, which make the final proofs complicated. Furthermore, Shamir's result [21] only works if a $1/\mathsf{poly}(d)$-accurate warm start is given, and the time needed to find a warm start is unclear.

In this paper, we develop a new algorithmic framework to solve $k$-SVD. It not only improves the aforementioned breakthroughs, but also relies only on simple convex analysis.

## 1.1 Our Results and the Settlement of an Open Question

We propose to use an extremely simple framework that we call LazySVD to solve $k$-SVD:

LazySVD: perform 1-SVD repeatedly, $k$ times in total.

More specifically, in this framework we first compute the leading singular vector $v$ of $A$, and then left-project $(I - vv^\top)A$ and repeat this procedure $k$ times. Quite surprisingly,

*This seemingly "most-intuitive" approach was widely considered as "not a good idea."*

In textbooks and research papers, one typically states that LazySVD has a running time that inversely depends on all the intermediate singular value gaps $\sigma_1 - \sigma_2, \dots, \sigma_k - \sigma_{k+1}$ [18, 21]. This dependence makes the algorithm useless if some singular values are close, and is even thought to be necessary [18]. For this reason, textbooks describe only block methods (such as block power method, block Krylov, alternating minimization) which find the top $k$ singular vectors together. Musco and Musco [19] stated as an *open question* to design "single-vector" methods without running time dependence on all the intermediate singular value gaps.

In this paper, we fully answer this open question with novel analyses on this LazySVD framework. In particular, the resulting running time either

- depends on $\mathsf{gap}^{-0.5}$ where $\mathsf{gap}$ is the relative singular value gap only between $\sigma_k$ and $\sigma_{k+1}$, or
- depends on $\varepsilon^{-0.5}$ where $\varepsilon$ is the approximation ratio (so is gap-free).

Such dependency matches the best known dependency for block methods.

More surprisingly, by making different choices of the 1-SVD subroutine in this LazySVD framework, we obtain multiple algorithms for different needs (see Table 1 and 2):

- If accelerated gradient descent or Lanczos algorithm is used for 1-SVD, we obtain a faster $k$-SVD algorithm than block Krylov [19].
- If a variance-reduction stochastic method is used for 1-SVD, we obtain the first accelerated stochastic algorithm for $k$-SVD, and this outperforms Shamir [21].
- If one sub-samples $A$ before applying LazySVD, the running time becomes $\widetilde{O}(\mathsf{nnz}(A) + \varepsilon^{-2}\mathsf{poly}(k) \cdot d)$. This improves upon [7] in certain (but sufficiently interesting) parameter regimes, but completely avoids the use of alternating minimization.

Finally, besides the running time advantages above, our analysis is completely based on convex optimization because 1-SVD is solvable using convex techniques. LazySVD also works when $k$ is not known to the algorithm, as opposed to block methods which need to know $k$ in advance.

| Paper | Running time | Frobenius norm | Spectral norm |
|---|---|---|---|
| [8] | $O(\mathsf{nnz}(A)) + O\big(\frac{k^2}{\varepsilon^4}(n+d) + \frac{k^3}{\varepsilon^5}\big)$ | $(1+\varepsilon)\|\Delta\|_F$ | $(1+\varepsilon)\|\Delta\|_F$ |
| [7] | $O(\mathsf{nnz}(A)) + \widetilde{O}\big(\frac{k^5(\sigma_1/\sigma_k)^2}{\varepsilon^2}(n+d)\big)$ | $(1+\varepsilon)\|\Delta\|_F$ | $\|\Delta\|_2 + \varepsilon\|\Delta\|_F$ |
| LazySVD **Theorem 5.1** | $O(\mathsf{nnz}(A)) + \widetilde{O}\big(\frac{k^2(\sigma_1/\sigma_{k+1})^4}{\varepsilon^2}d\big)$ | N/A | $\|\Delta\|_2 + \varepsilon\|\Delta\|_F$ |
| | $O(\mathsf{nnz}(A)) + \widetilde{O}\big(\frac{k^2(\sigma_1/\sigma_{k+1})^2}{\varepsilon^{2.5}}(n+d)\big)$ | N/A | $\|\Delta\|_2 + \varepsilon\|\Delta\|_F$ |
| | $O(\mathsf{nnz}(A)) + \widetilde{O}\big(\frac{k^4(\sigma_1/\sigma_{k+1})^{4.5}}{\varepsilon^2}d\big)$ | $(1+\varepsilon)\|\Delta\|_2$ | $\|\Delta\|_2 + \varepsilon\|\Delta\|_F$ |

Table 2: Performance comparison among $O(\mathsf{nnz}(A) + \mathsf{poly}(1/\varepsilon))$ type of algorithms. Remark: we have not tried hard to improve the dependency with respect to $k$ or $(\sigma_1/\sigma_{k+1})$. See Remark 5.2.

**Other Related Work.** Some authors focus on the streaming or online model of 1-SVD [4, 15, 17] or $k$-SVD [3]. These algorithms are slower than offline methods. Unlike $k$-SVD, accelerated stochastic methods were previously known for 1-SVD [12, 13]. After this paper is published, LazySVD has been generalized to also solve canonical component analysis and generalized PCA by the same authors [1]. If one is only interested in projecting a vector to the top $k$-eigenspace without computing the top $k$ eigenvectors like we do in this paper, this can also be done in an accelerated manner [2].

## 2 Preliminaries

Given a matrix $A$ we denote by $\|A\|_2$ and $\|A\|_F$ respectively the spectral and Frobenius norms of $A$. For $q \geq 1$, we denote by $\|A\|_{S_q}$ the Schatten $q$-norm of $A$. We write $A \succeq B$ if $A, B$ are symmetric and $A - B$ is positive semi-definite (PSD). We denote by $\lambda_k(M)$ the $k$-th largest eigenvalue of a symmetric matrix $M$, and $\sigma_k(A)$ the $k$-th largest singular value of a rectangular matrix $A$.

Since $\lambda_k(AA^\top) = \lambda_k(A^\top A) = (\sigma_k(A))^2$,

solving $k$-SVD for $A$ is the same as solving $k$-PCA for $M = AA^\top$.

We denote by $\sigma_1 \geq \cdots \sigma_d \geq 0$ the singular values of $A \in \mathbb{R}^{d\times n}$, by $\lambda_1 \geq \cdots \lambda_d \geq 0$ the eigenvalues of $M = AA^\top \in \mathbb{R}^{d\times d}$. (Although $A$ may have fewer than $d$ singular values for instance when $n < d$, if this happens, we append zeros.) We denote by $A_k^*$ the best rank-$k$ approximation of $A$.

We use $\perp$ to denote the orthogonal complement of a matrix. More specifically, given a column orthonormal matrix $U \in \mathbb{R}^{d\times k}$, we define $U^\perp \stackrel{\text{def}}{=} \{x \in \mathbb{R}^d \mid U^\top x = 0\}$. For notational simplicity, we sometimes also denote $U^\perp$ as a $d \times (d-k)$ matrix consisting of some basis of $U^\perp$.

**Theorem 2.1** (approximate matrix inverse)**.** *Given $d \times d$ matrix $M \succeq 0$ and constants $\lambda, \delta > 0$ satisfying $\lambda I - M \succeq \delta I$, one can minimize the quadratic $f(x) \stackrel{\text{def}}{=} x^\top(\lambda I - M)x - b^\top x$ in order to invert $(\lambda I - M)^{-1}b$. Suppose the desired accuracy is $\big\|x - (\lambda I - M)^{-1}b\big\| \leq \varepsilon$. Then,*

- *Accelerated gradient descent (AGD) produces such an output $x$ in $O\big(\frac{\lambda^{1/2}}{\delta^{1/2}}\log\frac{\lambda}{\varepsilon\delta}\big)$ iterations, each requiring $O(d)$ time plus the time needed to multiply $M$ with a vector.*
- *If $M$ is given in the form $M = \frac{1}{n}\sum_{i=1}^n a_i a_i^\top$ and $\|a_i\|_2 \leq 1$, then accelerated SVRG (see for instance [5]) produces such an output $x$ in time $O\big(\max\{nd, \frac{n^{3/4}d\lambda^{1/4}}{\delta^{1/2}}\}\log\frac{\lambda}{\varepsilon\delta}\big)$.*

**Algorithm 1** $\texttt{AppxPCA}(\mathcal{A}, M, \delta_\times, \varepsilon, p)$ $\diamond$ ***(only for proving our theoretical results; for practitioners, feel free to use your favorite 1-PCA algorithm such as Lanczos to replace*** `AppxPCA`***.)***

**Input:** $\mathcal{A}$, an approximate matrix inversion method; $M \in \mathbb{R}^{d\times d}$, a symmetric matrix satisfying $0 \preceq M \preceq I$; $\delta_\times \in (0, 0.5]$, a multiplicative error; $\varepsilon \in (0,1)$, a numerical accuracy parameter; and $p \in (0,1)$, a confidence parameter. $\diamond$ *running time only logarithmically depends on* $1/\varepsilon$ *and* $1/p$.

$m_1 \leftarrow \left\lceil 4\log\left(\frac{288d}{p^2}\right)\right\rceil$, $m_2 \leftarrow \left\lceil \log\left(\frac{36d}{p^2\varepsilon}\right)\right\rceil$;
$\diamond$ $m_1 = T^{\mathrm{PM}}(8, 1/32, p)$ *and* $m_2 = T^{\mathrm{PM}}(2, \varepsilon/4, p)$ *using definition in Lemma A.1*
$\widetilde{\varepsilon}_1 \leftarrow \frac{1}{64m_1}\left(\frac{\delta_\times}{6}\right)^{m_1}$ and $\widetilde{\varepsilon}_2 \leftarrow \frac{\varepsilon}{8m_2}\left(\frac{\delta_\times}{6}\right)^{m_2}$
$\widehat{w}_0 \leftarrow$ a random unit vector; $s \leftarrow 0$; $\lambda^{(0)} \leftarrow 1 + \delta_\times$;
**repeat**
    $s \leftarrow s + 1$;
    **for** $t = 1$ **to** $m_1$ **do**
        Apply $\mathcal{A}$ to find $\widehat{w}_t$ satisfying $\left\|\widehat{w}_t - (\lambda^{(s-1)}I - M)^{-1}\widehat{w}_{t-1}\right\| \leq \widetilde{\varepsilon}_1$;
    **end for**
    $w \leftarrow \widehat{w}_{m_1}/\|\widehat{w}_{m_1}\|$;
    Apply $\mathcal{A}$ to find $v$ satisfying $\left\|v - (\lambda^{(s-1)}I - M)^{-1}w\right\| \leq \widetilde{\varepsilon}_1$;
    $\Delta^{(s)} \leftarrow \frac{1}{2}\cdot\frac{1}{w^\top v - \widetilde{\varepsilon}_1}$ and $\lambda^{(s)} \leftarrow \lambda^{(s-1)} - \frac{\Delta^{(s)}}{2}$;
**until** $\Delta^{(s)} \leq \frac{\delta_\times \lambda^{(s)}}{3}$
$f \leftarrow s$;
**for** $t = 1$ **to** $m_2$ **do**
    Apply $\mathcal{A}$ to find $\widehat{w}_t$ satisfying $\left\|\widehat{w}_t - (\lambda^{(f)}I - M)^{-1}\widehat{w}_{t-1}\right\| \leq \widetilde{\varepsilon}_2$;
**end for**
**return** $w \stackrel{\text{def}}{=} \widehat{w}_{m_2}/\|\widehat{w}_{m_2}\|$.

# 3 A Specific 1-SVD Algorithm: Shift-and-Inverse Revisited

In this section, we study a specific 1-PCA algorithm `AppxPCA` (recall 1-PCA equals 1-SVD). It is a (multiplicative-)approximate algorithm for computing the leading eigenvector of a symmetric matrix.

We emphasize that, in principle, most known 1-PCA algorithms (e.g., power method, Lanczos method) are suitable for our LazySVD framework. We choose `AppxPCA` solely because it provides the maximum flexibility in obtaining all stochastic / NNZ running time results at once.

Our `AppxPCA` uses the shift-and-inverse routine [12, 13], and our pseudo-code in Algorithm 1 is a modification of Algorithm 5 that appeared in [12]. Since we need a more refined running time statement with a multiplicative error guarantee, and since the stated proof in [12] is anyways only a sketched one, we choose to carefully reprove a similar result of [12] (details in Appendix A) and state the following theorem:

**Theorem 3.1** (`AppxPCA`)**.** *Let* $M \in \mathbb{R}^{d\times d}$ *be a symmetric matrix with eigenvalues* $1 \geq \lambda_1 \geq \cdots \geq \lambda_d \geq 0$ *and corresponding eigenvectors* $u_1, \ldots, u_d$*. With probability at least* $1-p$*,* `AppxPCA` *produces an output* $w$ *satisfying*

$$\sum_{i\in[d],\lambda_i\leq(1-\delta_\times)\lambda_1} (w^\top u_i)^2 \leq \varepsilon \quad \textit{and} \quad w^\top M w \geq (1-\delta_\times)(1-\varepsilon)\lambda_1 \ .$$

*Furthermore, the total number of oracle calls to* $\mathcal{A}$ *is* $O(\log(1/\delta_\times)m_1 + m_2)$*, and each time we call* $\mathcal{A}$ *we have* $\frac{\lambda^{(s)}}{\lambda_{\min}(\lambda^{(s)}I-M)} \leq \frac{12}{\delta_\times}$ *and* $\frac{1}{\lambda_{\min}(\lambda^{(s)}I-M)} \leq \frac{12}{\delta_\times\lambda_1}$.

**Algorithm 2** $\texttt{LazySVD}(\mathcal{A}, M, k, \delta_\times, \varepsilon_{\mathsf{pca}}, p)$

**Input:** $\mathcal{A}$, an approximate matrix inversion method; $M \in \mathbb{R}^{d\times d}$, a matrix satisfying $0 \preceq M \preceq I$; $k \in [d]$, the desired rank; $\delta_\times \in (0,1)$, a multiplicative error; $\varepsilon_{\mathsf{pca}} \in (0,1)$, a numerical accuracy parameter; and $p \in (0,1)$, a confidence parameter.

$M_0 \leftarrow M$ and $V_0 \leftarrow []$;
**for** $s = 1$ **to** $k$ **do**
  $v'_s \leftarrow \texttt{AppxPCA}(\mathcal{A}, M_{s-1}, \delta_\times/2, \varepsilon_{\mathsf{pca}}, p/k)$;
    ⋄ ***to practitioners: use your favorite 1-PCA algorithm such as Lanczos to compute $v'_s$***
  $v_s \leftarrow \big((I - V_{s-1}V_{s-1}^\top)v'_s\big) / \big\|(I - V_{s-1}V_{s-1}^\top)v'_s\big\|$;   ⋄ *project $v'_s$ to $V_{s-1}^\perp$*
  $V_s \leftarrow [V_{s-1}, v_s]$;
  $M_s \leftarrow (I - v_s v_s^\top) M_{s-1} (I - v_s v_s^\top)$   ⋄ *we also have $M_s = (I - V_s V_s^\top) M (I - V_s V_s^\top)$*
**end for**
**return** $V_k$.

Since $\texttt{AppxPCA}$ reduces 1-PCA to oracle calls of a matrix inversion subroutine $\mathcal{A}$, the stated conditions $\frac{\lambda^{(s)}}{\lambda_{\min}(\lambda^{(s)}I - M)} \leq \frac{12}{\delta_\times}$ and $\frac{1}{\lambda_{\min}(\lambda^{(s)}I - M)} \leq \frac{12}{\delta_\times \lambda_1}$ in Theorem 3.1, together with complexity results for matrix inversions (see Theorem 2.1), imply the following running times for $\texttt{AppxPCA}$:

**Corollary 3.2.**

- *If $\mathcal{A}$ is AGD, the running time of* $\texttt{AppxPCA}$ *is $\widetilde{O}\big(\frac{1}{\delta_\times^{1/2}}\big)$ multiplied with $O(d)$ plus the time needed to multiply $M$ with a vector.*
- *If $M = \frac{1}{n}\sum_{i=1}^n a_i a_i^\top$ where each $\|a_i\|_2 \leq 1$, and $\mathcal{A}$ is accelerated SVRG, then the total running time of* $\texttt{AppxPCA}$ *is $\widetilde{O}\big(\max\{nd, \frac{n^{3/4} d}{\lambda_1^{1/4}\delta_\times^{1/2}}\}\big)$.*

# 4 Main Algorithm and Theorems

Our algorithm $\texttt{LazySVD}$ is stated in Algorithm 2. It starts with $M_0 = M$, and repeatedly applies $k$ times $\texttt{AppxPCA}$. In the $s$-th iteration, it computes an approximate leading eigenvector of matrix $M_{s-1}$ using $\texttt{AppxPCA}$ with a multiplicative error $\delta_\times/2$, projects $M_{s-1}$ to the orthogonal space of this vector, and then calls it matrix $M_s$.

In this stated form, $\texttt{LazySVD}$ finds approximately the top $k$ eigenvectors of a symmetric matrix $M \in \mathbb{R}^{d\times d}$. If $M$ is given as $M = AA^\top$, then $\texttt{LazySVD}$ automatically finds the $k$-SVD of $A$.

## 4.1 Our Core Theorems

We state our approximation and running time core theorems of $\texttt{LazySVD}$ below, and then provide corollaries to translate them into gap-dependent and gap-free theorems on $k$-SVD.

**Theorem 4.1** (approximation)**.** *Let $M \in \mathbb{R}^{d\times d}$ be a symmetric matrix with eigenvalues $1 \geq \lambda_1 \geq \cdots \lambda_d \geq 0$ and corresponding eigenvectors $u_1, \ldots, u_d$. Let $k \in [d]$, let $\delta_\times, p \in (0,1)$, and let $\varepsilon_{\mathsf{pca}} \leq \mathsf{poly}\big(\varepsilon, \delta_\times, \frac{1}{d}, \frac{\lambda_1}{\lambda_{k+1}}\big)$.[3] Then,* $\texttt{LazySVD}$ *outputs a (column) orthonormal matrix $V_k = (v_1, \ldots, v_k) \in \mathbb{R}^{d\times k}$ which, with probability at least $1 - p$, satisfies all of the following properties. (Denote by $M_k = (I - V_k V_k^\top) M (I - V_k V_k^\top)$.)*

[3]The detailed specifications of $\varepsilon_{\mathsf{pca}}$ can be found in the appendix where we restate the theorem more formally. To provide the simplest proof, we have not tightened the polynomial factors in the theoretical upper bound of $\varepsilon_{\mathsf{pca}}$ because the running time depends only logarithmic on $1/\varepsilon_{\mathsf{pca}}$.

(a) *Core lemma:* $\|V_k^\top U\|_2 \le \varepsilon$, *where* $U = (u_j, \dots, u_d)$ *is the (column) orthonormal matrix and* $j$ *is the smallest index satisfying* $\lambda_j \le (1-\delta_\times)\|M_{k-1}\|_2$.

(b) *Spectral norm guarantee:* $\lambda_{k+1} \le \|M_k\|_2 \le \frac{\lambda_{k+1}}{1-\delta_\times}$.

(c) *Rayleigh quotient guarantee:* $(1-\delta_\times)\lambda_k \le v_k^\top M v_k \le \frac{1}{1-\delta_\times}\lambda_k$.

(d) *Schatten-q norm guarantee: for every* $q \ge 1$, *we have* $\|M_k\|_{S_q} \le \frac{(1+\delta_\times)^2}{(1-\delta_\times)^2}\Big(\sum_{i=k+1}^d \lambda_i^q\Big)^{1/q}$.

We defer the proof of Theorem 4.1 to the appendix, but highlight the main ideas and techniques behind the proof in Section 4.3. Below we state the running time of `LazySVD`.

**Theorem 4.2** (running time)**.** `LazySVD` *can be implemented to run in time*

- $\widetilde{O}\big(\frac{k\mathsf{nnz}(M)+k^2 d}{\delta_\times^{1/2}}\big)$ *if* $\mathcal{A}$ *is AGD and* $M \in \mathbb{R}^{d\times d}$ *is given explicitly;*
- $\widetilde{O}\big(\frac{k\mathsf{nnz}(A)+k^2 d}{\delta_\times^{1/2}}\big)$ *if* $\mathcal{A}$ *is AGD and* $M$ *is given as* $M = AA^\top$ *where* $A \in \mathbb{R}^{d\times n}$*; or*
- $\widetilde{O}\big(knd + \frac{kn^{3/4}d}{\lambda_k^{1/4}\delta_\times^{1/2}}\big)$ *if* $\mathcal{A}$ *is accelerated SVRG and* $M = \frac{1}{n}\sum_{i=1}^n a_i a_i^\top$ *where each* $\|a_i\|_2 \le 1$.

*Above, the* $\widetilde{O}$ *notation hides logarithmic factors with respect to* $k, d, 1/\delta_\times, 1/p, 1/\lambda_1, \lambda_1/\lambda_k$.

*Proof of Theorem 4.2.* We call $k$ times `AppxPCA`, and each time we can feed $M_{s-1} = (I - V_{s-1}V_{s-1}^\top)M(I - V_{s-1}V_{s-1}^\top)$ implicitly into `AppxPCA` thus the time needed to multiply $M_{s-1}$ with a $d$-dimensional vector is $O(dk + \mathsf{nnz}(M))$ or $O(dk + \mathsf{nnz}(A))$. Here, the $O(dk)$ overhead is due to the projection of a vector into $V_{s-1}^\perp$. This proves the first two running times using Corollary 3.2.

To obtain the third running time, when we compute $M_s$ from $M_{s-1}$, we explicitly project $a_i' \leftarrow (I - v_s v_s^\top)a_i$ for each vector $a_i$, and feed the new $a_1', \dots, a_n'$ into `AppxPCA`. Now the running time follows from the second part of Corollary 3.2 together with the fact that $\|M_{s-1}\|_2 \ge \|M_{k-1}\|_2 \ge \lambda_k$. □

## 4.2 Our Main Results for $k$-SVD

Our main theorems imply the following corollaries (proved in Appendix C.1 for completeness).

**Corollary 4.3** (Gap-dependent $k$-SVD)**.** *Let* $A \in \mathbb{R}^{d\times n}$ *be a matrix with singular values* $1 \ge \sigma_1 \ge \cdots \sigma_d \ge 0$ *and the corresponding left singular vectors* $u_1, \dots, u_d \in \mathbb{R}^d$. *Let* $\mathsf{gap} = \frac{\sigma_k - \sigma_{k+1}}{\sigma_k}$ *be the relative gap. For fixed* $\varepsilon, p > 0$, *consider the output*

$$V_k \leftarrow \texttt{LazySVD}\Big(\mathcal{A}, AA^\top, k, \mathsf{gap}, O\big(\tfrac{\varepsilon^4\cdot \mathsf{gap}^2}{k^4(\sigma_1/\sigma_k)^4}\big), p\Big) \quad .$$

*Then, defining* $W = (u_{k+1}, \dots, u_d)$, *we have with probability at least* $1-p$:

$$V_k \text{ is a rank-}k\text{ (column) orthonormal matrix with } \quad \|V_k^\top W\|_2 \le \varepsilon \quad .$$

*Our running time is* $\widetilde{O}\big(\frac{k\mathsf{nnz}(A)+k^2 d}{\sqrt{\mathsf{gap}}}\big)$, *or time* $\widetilde{O}\big(knd + \frac{kn^{3/4}d}{\sigma_k^{1/2}\sqrt{\mathsf{gap}}}\big)$ *in the stochastic setting (1.1).*

Above, both running times depend only poly-logarithmically on $1/\varepsilon$.

**Corollary 4.4** (Gap-free $k$-SVD). *Let $A \in \mathbb{R}^{d\times n}$ be a matrix with singular values $1 \geq \sigma_1 \geq \cdots \sigma_d \geq 0$. For fixed $\varepsilon, p > 0$, consider the output*

$$(v_1, \dots, v_k) = V_k \leftarrow \texttt{LazySVD}\Big(\mathcal{A}, AA^\top, k, \tfrac{\varepsilon}{3}, O\Big(\tfrac{\varepsilon^6}{k^4 d^4 (\sigma_1/\sigma_{k+1})^{12}}\Big), p\Big) \enspace .$$

*Then, defining $A_k = V_k V_k^\top A$ which is a rank $k$ matrix, we have with probability at least $1-p$:*

1. *Spectral norm guarantee: $\|A - A_k\|_2 \leq (1+\varepsilon)\|A - A_k^*\|_2$;*
2. *Frobenius norm guarantee: $\|A - A_k\|_F \leq (1+\varepsilon)\|A - A_k^*\|_F$; and*
3. *Rayleigh quotient guarantee: $\forall i \in [k]$, $\big|v_i^\top AA^\top v_i - \sigma_i^2\big| \leq \varepsilon\sigma_i^2$.*

*Running time is $\widetilde{O}\big(\frac{k\mathsf{nnz}(A)+k^2 d}{\sqrt{\varepsilon}}\big)$, or time $\widetilde{O}\big(knd + \frac{kn^{3/4}d}{\sigma_k^{1/2}\sqrt{\varepsilon}}\big)$ in the stochastic setting (1.1).*

**Remark 4.5.** The spectral and Frobenius guarantees we adopted are standard. It was observed that the spectral guarantee is more desirable than the Frobenius one in practice [19]. In fact, our algorithm implies for all $q \geq 1$, $\|A - A_k\|_{S_q} \leq (1+\varepsilon)\|A - A_k^*\|_{S_q}$ where $\|\cdot\|_{S_q}$ is the Schatten-$q$ norm. Rayleigh-quotient types of guarantee were introduced by Musco and Musco [19] for a more refined comparison. They showed that the block Krylov method satisfies $|v_i^\top AA^\top v_i - \sigma_i^2| \leq \varepsilon\sigma_{k+1}^2$, which is slightly stronger than ours. However, these two guarantees are not much different in practice as we evidenced in our experiments.

## 4.3 High-Level Ideas Behind Our Theorems

For the sake of demonstrating the idea, we focus on the case when there is a (known) relative gap $\mathsf{gap} \stackrel{\text{def}}{=} (\sigma_k - \sigma_{k+1})/\sigma_k \in [0,1]$ between the $k$-th and the $(k+1)$-th singular values of $A$. Note that LazySVD can be equivalently viewed as follows. At iteration $s$, LazySVD starts with a (column) orthonormal matrix $V_{s-1} \in \mathbb{R}^{d\times(s-1)}$. It finds an approximate leading eigenvector $v_s$ of $M_{s-1} = (I - V_{s-1}V_{s-1}^\top)AA^\top(I - V_{s-1}V_{s-1}^\top)$, where $M_{s-1}$ is the projection of $M = AA^\top$ into space $V_{s-1}^\perp$. Then, LazySVD appends $V_s \leftarrow [V_{s-1}, v_s]$ and continues to the next iteration.

**Obtain Faster Running Time.** Ideally, if each $v_s$ were exactly the leading eigenvector of $M_s$, then our final $V_k$ would become exactly the top $k$ singular vectors of $A$. However, computing exact eigenvectors is too slow, so the main challenge is to tolerate as much error as possible to compute each $v_s$, in order to *tradeoff* for a faster running time.

It was previously a folklore that one can approximately compute each $v_s$ to a good precision so that the running time depends on *all* intermediate gaps $\frac{\sigma_i}{\sigma_i - \sigma_{i+1}}$ for $i = 1, \dots, k$. This is too slow, although thought to be somewhat necessary by some authors. A weaker alternative is to compute each $v_s$ to an additive precision so $v_s^\top M_{s-1} v_s \geq \sigma_s^2 - \mathsf{gap}\cdot\sigma_k^2$. However, this remains too slow because the running time would polynomially depend on $\sigma_1/\sigma_k$.

In LazySVD, we tolerate the $s$-th leading eigenvalue computation to suffer from a multiplicative error $\mathsf{gap}$ —more precisely, to satisfy $v_s^\top M_{s-1} v_s \geq (1-\mathsf{gap})\|M_{s-1}\|_2$. This implies our declared running time in Table 1 owing to the 1-PCA result of Section 3. What it remains is to prove the correctness of our algorithm: that is, to prove $\|M_{s-1}\|_2 \approx \sigma_s^2$ for each $s = 1, 2, \dots, k$.

**Obtain Correctness.** Our main idea is to use the fact that each vector $v_s$ "approximately" lies in the span of the top $k$ eigenvectors of $M = AA^\top$.

Notice if each $v_s$ perfectly lay in the span of the top $k$ eigenvectors of $M$, we would be able to claim —using the Cauchy interlacing theorem— that (1) the $(k+1-s)$-th largest eigenvalue of matrix $M_s$ would be *never be larger than* $\sigma_{k+1}^2$, and (2) the largest eigenvalue of $M_s$ would be never smaller than $\sigma_{s+1}^2$:

**Proposition 4.6** (Cauchy interlacing theorem)**.** *Given a symmetric matrix $N \in \mathbb{R}^{d\times d}$ with eigenvalues $\lambda_1 \geq \cdots \geq \lambda_d$, if $v$ is in the span of the first $k$ eigenvectors of $N$, then $(I - vv^\top)N(I - vv^\top)$ has eigenvalues $\lambda'_1 \geq \cdots \geq \lambda'_{d-1}$ satisfying:*

$$\forall i \in [k-1]\colon \lambda_i \geq \lambda'_i \geq \lambda_{i+1} \qquad \textit{and} \qquad \forall i \geq k\colon \lambda_{i+1} = \lambda'_i \enspace .$$

Therefore, the difference between the largest and the $(k+1-s)$-th largest eigenvalue of $M_s$ would be at least $\sigma_{s+1}^2 - \sigma_{k+1}^2 > \Omega(\mathsf{gap}) \cdot \sigma_{s+1}^2$ as long as $s = 0, 1, \ldots, k-1$. This could allow us to apply a "gap-multiplicative error" algorithm to obtain the next leading eigenvector $v_{s+1}$ of $M_s$, and this $v_{s+1}$ would "almost completely" lie in the span of the top $(k-s)$ eigenvectors of $M_s$ and thus the span of the top $k$ eigenvectors of $M$. Repeating this argument for $k$ times, we could have obtained the top $k$ singular vectors of $A$, up to rotation.

Our *main technique contribution* is to extend the above argument into an approximate setting, and to propagate error "moderately" across iterations. While a naive bound could easily blow up the error exponentially with respect to $k$, we provide non-trivial analysis to show that the error grows at most linearly in $k$. This step of our proof essentially consists of two parts.

Part 1: We show that, in each iteration, the vector $v_s$ obtained from `AppxPCA` only correlates with the last $d-k$ eigenvectors of $M_{s-1}$ by a polynomially small factor (i.e., by $\mathsf{poly}(\varepsilon, 1/n, 1/d)$).

Part 2: We develop a gap-free variant of the Wedin theorem for matrices [23], which translates Part 1 into two statements.

- The last $d-k$ eigenvectors of $M_s = (I - vv^\top)M_{s-1}(I - vv^\top)$ approximately lie in the span of the last $d-k$ eigenvectors of $M_{s-1}$.
- The $(k+1-s)$-th eigenvalue of $M_s$ is close to the $(k+2-s)$-th eigenvalue of $M_{s-1}$.

Recursively applying the above statements $s$ times, we conclude that

- The last $d-k$ eigenvectors of $M_s$ approximately lie in the span of the last $d-k$ eigenvectors of $M_0 = M$.
- The $(k+1-s)$-th eigenvalue of $M_s$ is close to $\sigma_{k+1}^2$, the $(k+1)$-th eigenvalue of $M_0 = M$.

These two properties replace the use of the Cauchy interlacing theorem, which only holds when "$v_s$ perfectly lay in the span of the top $k$ eigenvectors of $M$." They together imply that "the difference between the largest and the $(k+1-s)$-th largest eigenvalue of $M_s$" is at least $\gtrsim \sigma_{s+1}^2 - \sigma_{k+1}^2 > \Omega(\mathsf{gap}) \cdot \sigma_{s+1}^2$. Therefore, we can proceed to the next iteration $s+1$, and compute a gap-multiplicative approximate leading eigenvector $v_{s+1}$ of $M_s$.

This summarizes our intuition behind the correctness of `LazySVD`.

# 5 NNZ Running Time

In this section, we translate our results in the previous section into the $O(\mathsf{nnz}(A) + \mathsf{poly}(k, 1/\varepsilon)(n + d))$ running-time statements. The idea is surprisingly simple: we sample either random columns of $A$, or random entries of $A$, and then apply `LazySVD` to compute the $k$-SVD. Such translation directly gives either $1/\varepsilon^{2.5}$ results if AGD is used as the convex subroutine and either column or entry sampling is used, or a $1/\varepsilon^2$ result if accelerated SVRG and column sampling are used together.

We only informally state our theorem and defer all the details to Appendix D.

**Theorem 5.1** (informal)**.** *Let $A \in \mathbb{R}^{d\times n}$ be a matrix with singular values $\sigma_1 \geq \cdots \geq \sigma_d \geq 0$. For every $\varepsilon \in (0, 1/2)$, one can apply* `LazySVD` *with appropriately chosen $\delta_\times$ on a "carefully sub-sampled version" of $A$. Then, the resulting matrix $V \in \mathbb{R}^{d\times k}$ can satisfy*

- *spectral norm guarantee: $\|A - VV^\top A\|_2 \leq \|A - A_k^*\|_2 + \varepsilon\|A - A_k^*\|_F$;*[4]
- *Frobenius norm guarantee: $\|A - VV^\top A\|_F \leq (1 + \varepsilon)\|A - A_k^*\|_F$.*

*The total running time depends on (1) whether column or entry sampling is used, (2) which matrix inversion routine $\mathcal{A}$ is used, and (3) whether spectral or Frobenius guarantee is needed. We list our deduced results in Table 2 and the formal statements can be found in Theorem D.4, D.6, and D.9.*

**Remark 5.2.** The main purpose of our NNZ results is to demonstrate the strength of `LazySVD` framework in terms of improving the $\varepsilon$ dependency to $1/\varepsilon^2$. Since the $1/\varepsilon^2$ rate matches sampling complexity, it is *very challenging* have an NNZ result with $1/\varepsilon^2$ dependency.[5] We have not tried hard, and believe it possible, to improve the polynomial dependence with respect to $k$ or $(\sigma_1/\sigma_{k+1})$. Also, somewhat surprisingly, in our analysis Frobenius norms become harder to minimize as opposed to spectral norms; this is in contrast to known literatures where usually Frobenius results are easier to prove [7].

# 6 Experiments

We demonstrate the practicality of our LazySVD framework, and compare it to block power method or block Krylov method. We emphasize that in theory, the best worse-cast complexity for 1-PCA is obtained by `AppxPCA` on top of accelerated SVRG. However, for the size of our chosen datasets, Lanczos method runs faster than `AppxPCA` and therefore we adopt Lanczos method as the 1-PCA method for our LazySVD framework.[6]

**Datasets.** We use datasets `SNAP/amazon0302`, `SNAP/email-enron`, and `news20` that were also used by Musco and Musco [19], as well as an additional but famous dataset `RCV1`. The first two can be found on the SNAP website [16] and the last two can be found on the LibSVM website [11]. The four datasets give rise sparse matrices of dimensions $257570 \times 262111, 35600 \times 16507, 11269 \times 53975$, and $20242 \times 47236$ respectively.

**Implemented Algorithms.** For the block Krylov method, it is a well-known issue that the Lanczos type of three-term recurrence update is numerically unstable. This is why Musco and Musco [19] only used the stable variant of block Krylov which requires an orthogonalization of each $n \times k$ matrix with respect to all previously obtained $n \times k$ matrices. This greatly improves the numerical stability albeit sacrificing running time. We implement both these algorithms. In sum, we have implemented:

- PM: block power method for $T$ iterations.
- Krylov: stable block Krylov method for $T$ iterations [19].

[4] This is the best known spectral guarantee one can obtain using NNZ running time [7]. It is an open question whether the stricter $\|A - VV^\top A\|_2 \leq (1 + \varepsilon)\|A - A_k^*\|_2$ type of spectral guarantee is possible.

[5] On one hand, one can use dimension reduction such as [9] to reduce the problem size to $O(k/\varepsilon^2)$; to the best of our knowledge, it is impossible to obtain any NNZ result faster than $1/\varepsilon^3$ using solely dimension reduction. On the other hand, obtaining $1/\varepsilon^2$ dependency was the main contribution of [7]: they relied on alternating minimization but we have avoided it in our paper.

[6] Our LazySVD framework turns every 1-PCA method satisfying Theorem 3.1 (including Lanczos method) into a $k$-SVD solver. However, our theoretical results (esp. stochastic and NNZ) rely on `AppxPCA` because Lanczos is not a stochastic method.

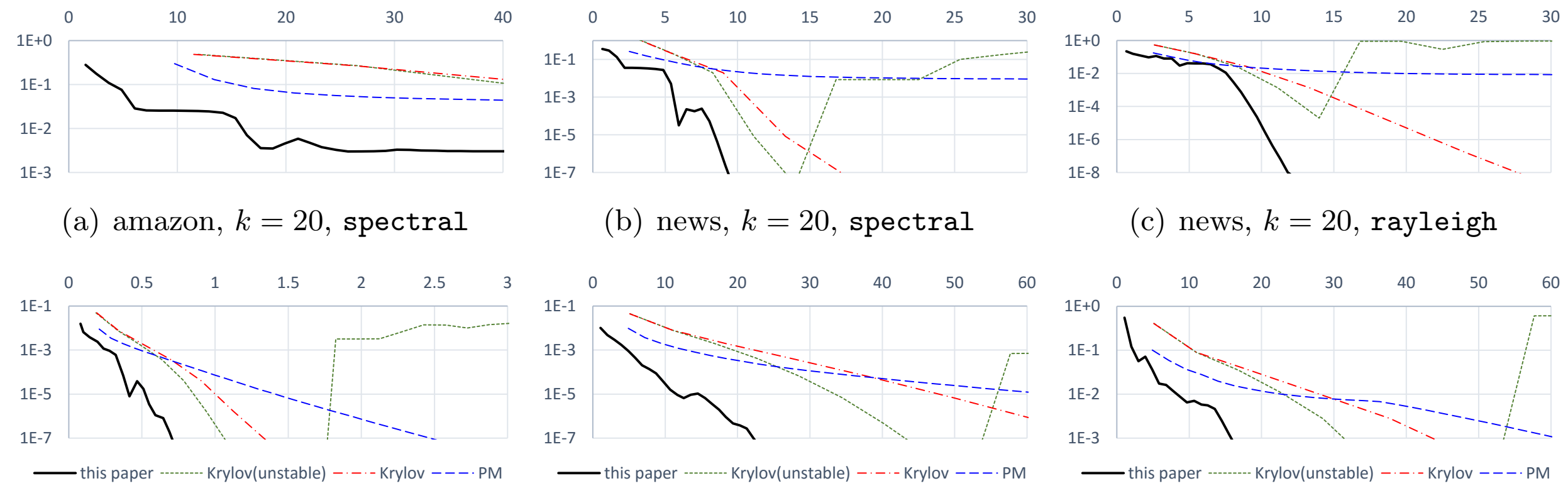

(a) amazon, $k = 20$, `spectral` (b) news, $k = 20$, `spectral` (c) news, $k = 20$, `rayleigh`

(d) email, $k = 10$, `Fnorm` (e) rcv1, $k = 30$, `Fnorm` (f) rcv1, $k = 30$, `rayleigh(last)`

Figure 1: Selected performance plots. Relative error ($y$-axis) vs. running time ($x$-axis).

- Krylov(unstable): the three-term recurrence implementation of block Krylov for $T$ iterations.
- `LazySVD`: $k$ calls of the vanilla Lanczos method, and each call runs $T$ iterations.

**A Fair Running-Time Comparison.** For a fixed integer $T$, the four methods go through the dataset (in terms of multiplying $A$ with column vectors) the same number of times. However, since `LazySVD` does not need block orthogonalization (as needed in PM and Krylov) and does not need a $(Tk)$-dimensional SVD computation in the end (as needed in Krylov), the running time of `LazySVD` is clearly much faster for a fixed value $T$. We therefore compare the performances of the four methods in terms of running time rather than $T$.

We programmed the four algorithms using the same programming language with the same sparse-matrix implementation. We tested them single-threaded on the same Intel i7-3770 3.40GHz personal computer. As for the final low-dimensional SVD decomposition step at the end of the PM or Krylov method (which is not needed for our `LazySVD`), we used a third-party library that is built upon the x64 Intel Math Kernel Library so the time needed for such SVD is maximally reduced.

**Performance Metrics.** We compute four metrics on the output $V = (v_1, \dots, v_k) \in \mathbb{R}^{n \times k}$:

- `Fnorm`: relative Frobenius norm error: $(\|A - VV^\top A\|_F - \|A - A_k^*\|_F)/\|A - A_k^*\|_F$.
- `spectral`: relative spectral norm error: $(\|A - VV^\top A\|_2 - \|A - A_k^*\|_2)/\|A - A_k^*\|_2$.
- `rayleigh(last)`: Rayleigh quotient error relative to $\sigma_{k+1}$: $\max_{j=1}^k \left|\sigma_j^2 - v_j^\top AA^\top v_j\right|/\sigma_{k+1}^2$.
- `rayleigh`: relative Rayleigh quotient error: $\max_{j=1}^k \left|\sigma_j^2 - v_j^\top AA^\top v_j\right|/\sigma_j^2$.

The first three metrics were also used by Musco and Musco [19]. We added the fourth one because our theory only predicted convergence with respect to the fourth but not the third metric. However, we observe that in practice they are not much different from each other.

**Our Results.** We study four datasets each with $k = 10, 20, 30$ and with the four performance metrics, totaling 48 plots. Due to space limitation, we only select six representative plots out of 48 and include them in Figure 1. (The full plots can be found in Figure 2, 3, 4 and 5 in the appendix.) We make the following observations:

- `LazySVD` outperforms its three competitors almost universally.
- Krylov(unstable) outperforms Krylov for small value $T$; however, it is less useful for obtaining accurate solutions due to its instability. (The dotted green curves even go up if $T$ is large.)

- Subspace power method performs the slowest unsurprisingly due to its lack of acceleration.

## Acknowledgements

We thank Ohad Shamir for listening to our ideas and pointing us to some of the existing results. The first author is partially supported by a Microsoft Research Award, no. 0518584, and an NSF grant, no. CCF-1412958.

# Appendix

## A Proof Details for Theorem 3.1: Convergence of AppxPCA

### A.1 Inexact Power Method

Consider the classical power method that is to start with a random unit vector $w_0$ and apply $w_t \leftarrow M w_{t-1}/\|M w_{t-1}\|$ iteratively.

**Lemma A.1** (Exact Power Method)**.** *Let $M$ be a PSD matrix with eigenvalues $\lambda_1 \geq \cdots \geq \lambda_d$ and the correpsonding eigenvectors $u_1, \dots, u_d$. Fix an error tolerance $\varepsilon > 0$, parameter $\kappa \geq 1$, and failure probability $p > 0$, define*

$$T^{\mathrm{PM}}(\kappa, \varepsilon, p) = \Big\lceil \frac{\kappa}{2} \log\Big(\frac{9d}{p^2 \varepsilon}\Big) \Big\rceil$$

*Then, with probability at least $1-p$ it holds that $\forall t \geq T^{\mathrm{PM}}(\kappa, \varepsilon, p)$:*

$$\sum_{i \in [d], \lambda_i \leq (1-1/\kappa)\lambda_1} (w_t^\top u_i)^2 \leq \varepsilon \quad \textit{and} \quad w_t^\top M w_t \geq (1 - 1/\kappa - \varepsilon)\lambda_1 \enspace .$$

*The probability of success depends only on the random variable $(w_0^\top u_1)^2$.*

*Proof.* For each $i \in [d]$,

$$\begin{aligned}(w_t^\top u_i)^2 &= \frac{\big((M^t w_0)^\top u_i\big)^2}{\|M^t w_0\|^2} = \frac{(w_0^\top M^t u_i)^2}{w_0^\top M^{2t} w_0} = \frac{\lambda_i^{2t} (w_0^\top u_i)^2}{\sum_{j=1}^d \lambda_j^{2t} (w_0^\top u_j)^2} = \frac{(w_0^\top u_i)^2}{\sum_{j=1}^d \big(\frac{\lambda_j}{\lambda_i}\big)^{2t} (w_0^\top u_j)^2} \\ &\leq \frac{(w_0^\top u_i)^2}{\big(\frac{\lambda_1}{\lambda_i}\big)^{2t} (w_0^\top u_1)^2} = \frac{(w_0^\top u_i)^2}{(w_0^\top u_1)^2} \Big(\frac{\lambda_i}{\lambda_1}\Big)^{2t}\end{aligned}$$

Since $w_0$ is a random unit vector, according to for instance Lemma 5 of [6], it holds with probability at least $1-p$ that $(w_0^\top u_1)^2 \geq \frac{p^2}{9d}$. Substituting this into the above inequality, we conclude that with probability at least $1-p$, for all $i \in [d]$, we have

$$(w_t^\top u_i)^2 \leq (w_0^\top u_i) \cdot \frac{9d}{p^2} \Big(1 - \frac{\lambda_1 - \lambda_i}{\lambda_1}\Big)^{2t} \leq (w_0^\top u_i) \cdot \frac{9d}{p^2} \cdot \exp\Big(-2 \frac{\lambda_1 - \lambda_i}{\lambda_1} t\Big)$$

As a result, for every $t \geq T^{\mathrm{PM}}(\kappa, \varepsilon, p)$, and every $i$ such that $\lambda_i \leq (1 - 1/\kappa)\lambda_1$ (which implies $\frac{\lambda_1 - \lambda_i}{\lambda_1} \geq 1/\kappa$), we have

$$(w_t^\top u_i)^2 \leq \varepsilon \cdot (w_0^\top u_i)^2$$

Summing them up we have

$$\sum_{i\in[d],\lambda_i\le(1-1/\kappa)\lambda_1}(w_t^\top u_i)^2\le\varepsilon\sum_{i\in[d]}(w_0^\top u_i)^2=\varepsilon\ .$$

This finishes the proof of the first bound. To prove the second bound, we compute that

$$\begin{aligned}w_t^\top M w_t=\sum_{i=1}^d\lambda_i(w_t^\top u_i)^2&\ge\sum_{i\in[d],\lambda_i>(1-1/\kappa)\lambda_1}\lambda_i(w_t^\top u_i)^2\ge(1-1/\kappa)\lambda_1\cdot\sum_{i\in[d],\lambda_i>(1-1/\kappa)\lambda_1}(w_t^\top u_i)^2\\&\ge(1-1/\kappa)(1-\varepsilon)\lambda_1\ge(1-1/\kappa-\varepsilon)\lambda_1\ .\end{aligned}$$

□

**Lemma A.2** (Lemma 4.1 of [12]). *Let $M$ be a PSD matrix with eigenvalues $\lambda_1\ge\cdots\lambda_d$. Fix an accuracy parameter $\widetilde{\varepsilon}>0$, and consider two update sequences*

$$\begin{aligned}&\widehat{w}_0^*=w_0,\qquad\forall t\ge1\colon\ \widehat{w}_t^*\leftarrow M\widehat{w}_{t-1}^*\\&\widehat{w}_0=w_0,\qquad\forall t\ge1\colon\ \widehat{w}_t\ \text{satisfies}\ \|\widehat{w}_t-M\widehat{w}_{t-1}\|\le\widetilde{\varepsilon},\end{aligned}$$

*Then, defining $w_t=\widehat{w}_t/\|\widehat{w}_t\|$ and $w_t^*=\widehat{w}_t^*/\|\widehat{w}_t^*\|$, it satisfies*

$$\|w_t-w_t^*\|\le\widetilde{\varepsilon}\cdot\Gamma(M,t),$$

*where*

$$\Gamma(M,t)\stackrel{\text{def}}{=}\frac{2}{\lambda_d^t}\begin{cases}t, & \text{if }\lambda_1=1;\\(\lambda_1^t-1)/(\lambda_1-1), & \text{if }\lambda_1\neq1.\end{cases}\quad\text{and we have }\Gamma(M,t)\le 2t\cdot\frac{\max\{1,\lambda_1^t\}}{\lambda_d^t}$$

**Theorem A.3** (Inexact Power Method). *Let $M$ be a PSD matrix with eigenvalues $\lambda_1\ge\cdots\ge\lambda_d$ and the corresponding eigenvectors $u_1,\dots,u_d$. With probability at least $1-p$ it holds that, for every $\varepsilon\in(0,1)$ and every $t\ge T^{\mathrm{PM}}(\kappa,\varepsilon/4,p)$, if $w_t$ is generated by the power method with per-iteration error $\widetilde{\varepsilon}=\frac{\varepsilon}{4\Gamma(M,t)}$, then*

$$\sum_{i\in[d],\lambda_i\le(1-1/\kappa)\lambda_1}(w_t^\top u_i)^2\le\varepsilon\quad\text{and}\quad w_t^\top M w_t\ge(1-1/\kappa-\varepsilon)\lambda_1\ .$$

*Proof.* Denoting by $w_t^*$ the output of power method with exact updates (with the same starting vector $w_0$ following Lemma A.2)

$$\begin{aligned}\sum_{j\in[d],\lambda_j\le(1-1/\kappa)\lambda_1}(w_t^\top u_j)^2&=\sum_{j\in[d],\lambda_j\le(1-1/\kappa)\lambda_1}(\langle w_t^*,u_j\rangle+\langle w_t-w_t^*,u_j\rangle)^2\\&\le\sum_{j\in[d],\lambda_j\le(1-1/\kappa)\lambda_1}2(\langle w_t^*,u_j\rangle)^2+2(\langle w_t-w_t^*,u_j\rangle)^2\\&\le\frac{\varepsilon}{2}+2\sum_{j\in[d]}(w_t-w_t^*)^\top u_ju_j^\top(w_t-w_t^*)=\frac{\varepsilon}{2}+2\|w_t-w_t^*\|^2\le\varepsilon\ .\end{aligned}$$

Above, the first inequality is because $(a+b)^2\le2a^2+2b^2$, the second inequality is due to the definition of $w_t^*$ and Lemma A.1, and the last inequality is because $\|w_t-w_t^*\|\le\widetilde{\varepsilon}\cdot\Gamma(M,t)=\frac{\varepsilon}{4}$ which implies $2\|w_t-w_t^*\|^2\le\varepsilon^2/8<\varepsilon/2$.

This finishes the proof of the first bound. The proof of the second bound is identical to the last paragraph of the proof of Lemma A.1. □

### A.2 Proof of Theorem 3.1

**Lemma A.4.** *With probability at least $1-p$, it holds that (where $\lambda_1$ is the largest eigenvalue of $M$.)*

*(a)* $\widetilde{\varepsilon}_1 \le \frac{1}{32\Gamma((\lambda^{(s-1)}I-M)^{-1},m_1)}$ *for each iteration $s$;*

*(b)* $\widetilde{\varepsilon}_2 \le \frac{\varepsilon}{4\Gamma((\lambda^{(f)}I-M)^{-1},m_2)}$*;*

*(c)* $0 \le \frac{1}{2}(\lambda^{(s-1)}-\lambda_1) \le \Delta^{(s)} \le \lambda^{(s-1)}-\lambda_1$ *and* $\frac{1}{2}(\lambda^{(s-1)}-\lambda_1) \le \lambda^{(s)}-\lambda_1$ *for each iteration $s$; and*

*(d)* $\lambda^{(f)}-\lambda_1 \in [\frac{\delta_\times}{12}\lambda^{(f)}, \delta_\times\lambda_1]$ *when the repeat-until loop is over.*

*Proof.* We denote by $A^{(s)} \stackrel{\text{def}}{=} (\lambda^{(s)}I-M)^{-1}$ for notational simplicity. Below we prove all the items by induction for a specific iteration $s \ge 2$ assuming that the items of the previous $s-1$ iterations are true. The base case of $s=1$ can be verified similar to the general arguments below but requiring some non-trivial notational changes. We omitted the proof of the base case $s=1$ in this paper.

(a) Recall that

$$\Gamma(A^{(s-1)},t) \le 2t \cdot \frac{\max\{1, \lambda_{\max}(A^{(s-1)})^t\}}{\lambda_{\min}(A^{(s-1)})^t}$$

On one hand, we have $\lambda_{\max}(A^{(s-1)}) = \frac{1}{\lambda^{(s-1)}-\lambda_1} \le \frac{2}{\lambda^{(s-2)}-\lambda_1} \le \frac{2}{\Delta^{(s-1)}}$ using Lemma A.4.c of the previous iteration. Combining this with the termination criterion $\Delta^{(s-1)} \ge \frac{\delta_\times}{3}\lambda^{(s-1)}$, we have $\lambda_{\max}(A^{(s-1)}) \le \frac{6}{\delta_\times \lambda^{(s-1)}}$. On the other hand, we have $\lambda_{\min}(A^{(s-1)}) = \frac{1}{\lambda^{(s-1)}-\lambda_d} \ge \frac{1}{\lambda^{(s-1)}}$. Combining the two bounds we conclude that $\Gamma(A^{(s-1)},t) \le 2t(6/\delta_\times)^t$. It is now obvious that $\widetilde{\varepsilon}_1 \le \frac{1}{32\Gamma(A^{(s-1)},m_1)}$ is satisfied because $\widetilde{\varepsilon}_1 = \frac{1}{64m_1}\big(\frac{\delta_\times}{6}\big)^{m_1}$.

(b) The same analysis as in the proof of Lemma A.4.a suggests that $\Gamma(A^{(f)},t) \le 2t(6/\delta_\times)^t$. This immediately yields $\widetilde{\varepsilon}_2 \le \frac{\varepsilon}{4\Gamma(A^{(f)},m_2)}$ because $\widetilde{\varepsilon}_2 = \frac{\varepsilon}{8m_2}\big(\frac{\delta_\times}{6}\big)^{m_2}$.

(c) Because Lemma A.4.a holds for the current iteration $s$ we can apply Theorem A.3 (with $\varepsilon = 1/8$ and $\kappa = 8$) and get

$$w^\top A^{(s-1)} w \ge \frac{3}{4}\lambda_{\max}(A^{(s-1)}) \ .$$

By the definition of $v$ in `AppxPCA` and the Cauchy-Schwartz inequality it holds that

$$w^\top v = w^\top A^{(s-1)} w + w^\top\big(v - A^{(s-1)}w\big) \in \big[w^\top A^{(s-1)}w - \widetilde{\varepsilon}_1, w^\top A^{(s-1)}w + \widetilde{\varepsilon}_1\big]$$

Combining the above two equations we have

$$\begin{aligned} w^\top v - \widetilde{\varepsilon}_1 &\in \Big[\frac{3}{4}\lambda_{\max}(A^{(s-1)}) - 2\widetilde{\varepsilon}_1, \lambda_{\max}(A^{(s-1)})\Big] \\ &\subseteq \Big[\frac{1}{2}\lambda_{\max}(A^{(s-1)}), \lambda_{\max}(A^{(s-1)})\Big] = \Big[\frac{1}{2}, 1\Big] \cdot \frac{1}{\lambda^{(s-1)}-\lambda_1} \ . \end{aligned}$$

In other words, $\Delta^{(s)} \stackrel{\text{def}}{=} \frac{1}{2} \cdot \frac{1}{w^\top v - \widetilde{\varepsilon}_1} \in \big[\frac{1}{2}(\lambda^{(s-1)}-\lambda_1), \lambda^{(s-1)}-\lambda_1\big]$.

At the same time, our update rule $\lambda^{(s)} = \lambda^{(s-1)} - \Delta^{(s)}/2$ ensures that $\lambda^{(s)}-\lambda_1 = \lambda^{(s-1)} - \lambda_1 - \Delta^{(s)}/2 \ge \lambda^{(s-1)} - \lambda_1 - \frac{\lambda^{(s-1)}-\lambda_1}{2} = \frac{1}{2}(\lambda^{(s-1)}-\lambda_1)$.

(d) The upper bound holds because $\lambda^{(f)} - \lambda_1 = \lambda^{(f-1)} - \frac{\Delta^{(f)}}{2} - \lambda_1 \leq \frac{3}{2}\Delta^{(f)} \leq \frac{\delta_\times \lambda^{(f)}}{2}$ where the first inequality follows from Lemma A.4.c of this last iteration, and the second inequality follows from our termination criterion $\Delta^{(f)} \leq \frac{\delta_\times \lambda^{(f)}}{3}$. Simply rewriting this inequality we have $\lambda^{(f)} - \lambda_1 \leq \frac{\delta_\times/2}{1-\delta_\times/2}\lambda_1 \leq \delta_\times \lambda_1$.

The lower bound is because using Lemma A.4.c (of this and the previous iteration) we have $\lambda^{(f)} - \lambda_1 \geq \frac{1}{4}\big(\lambda^{(f-2)} - \lambda_1\big) \geq \frac{\Delta^{(f-1)}}{4} \overset{①}{\geq} \frac{\delta_\times \lambda^{(f-1)}}{12} \geq \frac{\delta_\times \lambda^{(f)}}{12}$. Here, inequality ① is because $\Delta^{(f-1)} > \frac{\delta_\times \lambda^{(f-1)}}{3}$ due to the termination criterion.

Finally since the success of Theorem A.3 only depends on the randomness of $\widehat{w}_0$, we have that with probability at least $1-p$ that all the above items are satisfied. □

*Proof of Theorem 3.1.* It follows from Theorem A.3 (with $\kappa = 2$) that, letting $\mu_i = 1/(\lambda^{(f)} - \lambda_i)$ be the $i$-th largest eigenvalue of the matrix $(\lambda^{(f)} I - M)^{-1}$, then

$$\sum_{i\in[d],\mu_i\leq\mu_1/2} (w^\top u_i)^2 \leq \varepsilon$$

Note that if an index $i \in [d]$ satisfies $\lambda_1 - \lambda_i \geq \delta_\times \lambda_1$, then we must have $\lambda_1 - \lambda_i \geq \lambda^{(f)} - \lambda_1$ owing to $\lambda^{(f)} - \lambda_1 \leq \delta_\times \lambda_1$ from Lemma A.4.d. This further implies that $2(\lambda^{(f)} - \lambda_1) \leq \lambda^{(f)} - \lambda_i$ and therefore $\mu_1/2 \geq \mu_i$. In sum, we also have

$$\sum_{i\in[d],\lambda_i\leq(1-\delta_\times)\lambda_1} (w^\top u_i)^2 \leq \varepsilon \ .$$

On the other hand,

$$w^\top M w = \sum_{i=1}^{d} \lambda_i (w^\top u_i)^2 \geq \sum_{i\in[d],\lambda_i>(1-\delta_\times)\lambda_1} \lambda_i (w^\top u_i)^2 \geq (1-\delta_\times)\lambda_1 \cdot \sum_{i\in[d],\lambda_i>(1-\delta_\times)\lambda_1} (w^\top u_i)^2$$
$$\geq (1-\delta_\times)(1-\varepsilon)\lambda_1 \ .$$

The number of oracle calls to $\mathcal{A}$ is determined by the number of iterations in the repeat-until loop. It is easy to verify that there are at most $O(\log(1/\delta_\times))$ such iteartions, so the total number of oracle calls to $\mathcal{A}$ is only $O(\log(1/\delta_\times) m_1 + m_2)$.

In addition, each time we call $\mathcal{A}$ we have

$$\frac{\lambda^{(s)}}{\lambda_{\min}(\lambda^{(s)} I - M)} \leq \frac{\lambda^{(s)}}{\lambda^{(s)} - \lambda_1}$$

If $s = 0$ then we have $\frac{\lambda^{(0)}}{\lambda^{(0)} - \lambda_1} \leq \frac{1+\delta_\times}{\delta_\times}$ because $\lambda_1 \leq 1$. If $s \leq f-2$ then we have $\frac{\lambda^{(s)}}{\lambda^{(s)} - \lambda_1} \leq \frac{\lambda^{(s)}}{\Delta^{(s+1)}} \leq \frac{\lambda^{(s)}}{\delta_\times \lambda^{(s+1)}/3} \leq \frac{3}{\delta_\times}$ where the first inequality follows from Lemma A.4.c, the second one follows from the stopping criterion, and the last one follows from the monotonicity of $\lambda^{(s)}$. If $s = f-1$ then we have $\frac{\lambda^{(s)}}{\lambda^{(s)} - \lambda_1} \leq \frac{2\lambda^{(s)}}{\lambda^{(s-1)} - \lambda_1} \leq \frac{2\lambda^{(s)}}{\Delta^{(s)}} \leq \frac{2\lambda^{(s)}}{\delta_\times \lambda^{(s)}/3} = \frac{6}{\delta_\times}$ where the first two inequalities follow from Lemma A.4.c and the third inequality follows from our stopping criterion. If $s = f$ then we have $\frac{\lambda^{(s)}}{\lambda^{(s)} - \lambda_1} \leq \frac{12}{\delta_\times}$ owing to Lemma A.4.d. In all cases we have $\frac{\lambda^{(s)}}{\lambda_{\min}(\lambda^{(s)} I - M)} \leq \frac{12}{\delta_\times}$.

Finally, we have $\frac{1}{\lambda_{\min}(\lambda^{(s)} I - M)} = \frac{\lambda^{(s)}}{\lambda_{\min}(\lambda^{(s)} I - M)} \cdot \frac{1}{\lambda^{(s)}} \leq \frac{12}{\delta_\times \lambda_1}$ where the last inequality follows from $\lambda^{(s)} \geq \lambda_1$. □

# B Lemmas Needed for Proving Our Main Theorem

In this section we provide some necessary lemmas on matrices that shall become essential for our proof of Theorem 4.1.

**Proposition B.1.** *Let $A, B$ be two (column) orthonormal matrix such that for $\eta \geq 0$,*

$$A^\top BB^\top A \succeq (1-\eta)I$$

*Then we have: there exists a matrix $Q$, $\|Q\|_2 \leq 1$ such that*

$$\|A - BQ\|_2 \leq \sqrt{\eta}$$

*Proof.* Since $A^\top A = I$ and $A^\top BB^\top A \succeq (1-\eta)I$, we know that $A^\top B^\perp (B^\perp)^\top A \preceq \eta I$. By the fact that

$$A = (BB^\top + B^\perp (B^\perp)^\top)A = BB^\top A + B^\perp (B^\perp)^\top A$$

we can let $Q = B^\top A$ and obtain

$$\|A - BQ\|_2 \leq \|B^\perp (B^\perp)^\top A\|_2 \leq \sqrt{\eta} \ . \qquad \square$$

## B.1 Approximate Projection Lemma

**Lemma B.2.** *Let $M$ be a PSD matrix with eigenvalues $\lambda_1 \geq \cdots \geq \lambda_d$ and the corresponding eigenvectors $u_1, \ldots, u_d \in \mathbb{R}^d$. For every $k \geq 1$, define $U^\perp = (u_1, \ldots, u_k) \in \mathbb{R}^{d\times k}$ and $U = (u_{k+1}, \ldots, u_d) \in \mathbb{R}^{d\times(d-k)}$. For every $\varepsilon \in (0, \frac{1}{2})$, let $V_s \in \mathbb{R}^{d\times s}$ be an orthogonal matrix such that $\|V_s^\top U\|_2 \leq \varepsilon$, define $Q_s \in \mathbb{R}^{d\times s}$ to be an arbitrary orthogonal basis of the column span of $U^\perp (U^\perp)^\top V_s$, then we have:*

$$\left\| \left(I - Q_s Q_s^\top\right) M \left(I - Q_s Q_s^\top\right) - \left(I - V_s V_s^\top\right) M \left(I - V_s V_s^\top\right) \right\|_2 \leq 13\lambda_1 \varepsilon \ .$$

*Proof of Lemma B.2.* Since $Q_s$ is an orthogonal basis of the column span of $U^\perp (U^\perp)^\top V_s$, there is a matrix $R \in \mathbb{R}^{s\times s}$ such that

$$Q_s = U^\perp (U^\perp)^\top V_s R$$

Using the fact that $Q_s^\top Q_s = I$, we have:

$$(U^\perp (U^\perp)^\top V_s R)^\top (U^\perp (U^\perp)^\top V_s R) = I \implies R^\top V_s^\top U^\perp (U^\perp)^\top V_s R = I \ .$$

By the fact that $V_s^\top V_s = I$ and $U^\perp (U^\perp)^\top + UU^\top = I$, we can rewrite the above equality as:

$$R^\top \left(I - V_s^\top UU^\top V_s\right) R = I \tag{B.1}$$

From our lemma assumption, we have: $\|V_s^\top U\|_2 \leq \varepsilon$, which implies $0 \preceq V_s^\top UU^\top V_s \preceq \varepsilon^2 I$. Putting this into (B.1), we obtain:

$$I \preceq R^\top R \preceq \frac{1}{1-\varepsilon^2} I \preceq \left(1 + \frac{4}{3}\varepsilon^2\right) I$$

The above inequality directly implies that $I \preceq RR^\top \preceq \big(1 + \frac{4}{3}\varepsilon^2\big)I$. Therefore,

$$
\begin{aligned}
& \left\|Q_s Q_s^\top - V_s V_s^\top\right\|_2 \\
= \;& \left\|U^\perp (U^\perp)^\top V_s R R^\top V_s^\top U^\perp (U^\perp)^\top - V_s V_s^\top\right\|_2 \\
= \;& \left\|U^\perp (U^\perp)^\top V_s R R^\top V_s^\top U^\perp (U^\perp)^\top - (U^\perp (U^\perp)^\top + UU^\top) V_s V_s^\top (U^\perp (U^\perp)^\top + UU^\top)\right\|_2 \\
\leq \;& \left\|U^\perp (U^\perp)^\top V_s (RR^\top - I) V_s^\top U^\perp (U^\perp)^\top\right\|_2 + \left\|UU^\top V_s V_s^\top UU^\top\right\|_2 + 2\left\|U^\perp (U^\perp)^\top V_s V_s^\top UU^\top\right\|_2 \\
\leq \;& \left\|RR^\top - I\right\|_2 + \left\|U^\top V_s V_s^\top U\right\|_2 + 2\left\|V_s^\top UU^\top V_s\right\|_2^{1/2} \\
\leq \;& \frac{4}{3}\varepsilon^2 + \varepsilon^2 + 2\varepsilon < \frac{19}{6}\varepsilon \ .
\end{aligned}
$$

Finally, we have

$$
\begin{aligned}
& \left\|\left(I - Q_s Q_s^\top\right) M \left(I - Q_s Q_s^\top\right) - \left(I - V_s V_s^\top\right) M \left(I - V_s V_s^\top\right)\right\|_2 \\
& \leq 2\left\|\left(Q_s Q_s^\top - V_s V_s^\top\right) M\right\|_2 + \left\|\left(Q_s Q_s^\top - V_s V_s^\top\right) M Q_s Q_s^\top\right\|_2 + \left\|\left(Q_s Q_s^\top - V_s V_s^\top\right) M V_s V_s^\top\right\|_2 \\
& \leq \frac{19 \times 4}{6}\lambda_1 \varepsilon < 13\lambda_1\varepsilon \ .
\end{aligned}
$$

□

## B.2 Gap-Free Wedin Theorem

**Lemma B.3** (Gap free Wedin Theorem)**.** *For $\varepsilon \geq 0$, let $A, B$ be two PSD matrices such that $\|A - B\|_2 \leq \varepsilon$. For every $\mu \geq 0$, $\tau > 0$, let $U$ be column orthonormal matrix consisting of eigenvectors of $A$ with eigenvalue $\leq \mu$, let $V$ be column orthonormal matrix consisting of eigenvectors of $B$ with eigenvalue $\geq \mu + \tau$, then we have:*

$$\|U^\top V\| \leq \frac{\varepsilon}{\tau} \ .$$

*Proof of Lemma B.3.* We write $A$ and $B$ in terms of eigenvalue decomposition:

$$A = U\Sigma U^\top + U'\Sigma' U'^\top \quad \text{and} \quad B = V\widetilde{\Sigma}V^\top + V'\widetilde{\Sigma}'V'^\top \ ,$$

where $U'$ is orthogonal to $U$ and $V'$ is orthogonal to $V$. Letting $R = A - B$, we obtain:

$$
\begin{aligned}
& \Sigma U^\top = U^\top A = U^\top (B + R) \\
\Longrightarrow\; & \Sigma U^\top V = U^\top B V + U^\top R V = U^\top V \widetilde{\Sigma} + U^\top R V \\
\Longrightarrow\; & \Sigma U^\top V \widetilde{\Sigma}^{-1} = U^\top V + U^\top R V \widetilde{\Sigma}^{-1} \ .
\end{aligned}
$$

Taking spectral norm on both sides, we obtain:

$$\|\Sigma\|_2 \|U^\top V\|_2 \|\widetilde{\Sigma}^{-1}\|_2 \geq \|\Sigma U^\top V \widetilde{\Sigma}^{-1}\|_2 \geq \|U^\top V\|_2 - \|U^\top R V \widetilde{\Sigma}^{-1}\|_2 \ .$$

This can be simplified to

$$\frac{\mu}{\mu + \tau}\|U^\top V\|_2 \geq \|U^\top V\|_2 - \frac{\varepsilon}{\mu + \tau} \ ,$$

and therefore we have $\|U^\top V\|_2 \leq \frac{\varepsilon}{\tau}$ as desired. □

### B.3 Projected Power Method

**Lemma B.4.** *Let $M \in \mathbb{R}^{d\times d}$ be a PSD matrix with eigenvalues $\lambda_1 \geq \cdots \geq \lambda_d$ and corresponding eigenvectors $u_1, \ldots, u_d$. Define $U = (u_{j+1}, \ldots, u_d) \in \mathbb{R}^{d\times(d-j)}$ to be the matrix consisting of all eigenvectors with eigenvalue $\leq \mu$. Let $v \in \mathbb{R}^d$ be a unit vector such that $\|v^\top U\|_2 \leq \varepsilon \leq 1/2$, and define*

$$M' = \left(I - vv^\top\right) M \left(I - vv^\top\right)$$

*Then, denoting $[V_2, V_1, v] \in \mathbb{R}^{d\times d}$ as the unitary matrix consisting of (column) eigenvectors of $M'$ with descending eigenvalues, where $V_1$ consists of eigenvectors with eigenvalue $\leq \mu + \tau$, then there exists a matrix $Q$, $\|Q\|_2 \leq 1$ such that*

$$\|U - V_1 Q\|_2 \leq \sqrt{\frac{169\lambda_1^2\varepsilon^2}{\tau^2} + \varepsilon^2}$$

*Proof of Lemma B.4.* Using Lemma B.2, let $q = \frac{U^\perp (U^\perp)^\top v}{\|U^\perp (U^\perp)^\top v\|_2}$ be the projection of $v$ to $U^\perp$, we know that

$$\left\|\left(I - qq^\top\right) M \left(I - qq^\top\right) - \left(I - vv^\top\right) M \left(I - vv^\top\right)\right\|_2 \leq 13\lambda_1\varepsilon \ .$$

Denote $\left(I - qq^\top\right) M \left(I - qq^\top\right)$ as $M''$. We know that $u_{j+1}, \ldots, u_d$ are still eigenvectors of $M''$ with eigenvalue $\lambda_{j+1}, \ldots, \lambda_d$.

Apply Lemma B.3 on $A = M''$, $U$ and $B = M'$, $V = V_2$, we obtain:

$$\|U^\top V_2\|_2 \leq \frac{13\lambda_1\varepsilon}{\tau} \ .$$

This implies that

$$U^\top V_1 V_1^\top U = I - U^\top V_2 V_2^\top U - U^\top vv^\top U \succeq \left(1 - \frac{169\lambda_1^2\varepsilon^2}{\tau^2} - \varepsilon^2\right) I \ ,$$

where the inequality uses the assumption $\|v^\top U\|_2 \leq \varepsilon$.

Apply Proposition B.1 to $A = U$ and $B = V_1$, we conclude that there exists a matrix $Q$, $\|Q\|_2 \leq 1$ such that

$$\|U - V_1 Q\|_2 \leq \sqrt{\frac{169\lambda_1^2\varepsilon^2}{\tau^2} + \varepsilon^2} \ .$$

□

## C Proof Details for Theorem 4.1: Our Main Theorem

In this section we prove Theorem 4.1 formally.

**Theorem 4.1** (restated). *Let $M \in \mathbb{R}^{d\times d}$ be a symmetric matrix with eigenvalues $1 \geq \lambda_1 \geq \cdots \lambda_d \geq 0$ and corresponding eigenvectors $u_1, \dots, u_d$. Let $k \in [d]$, and $\delta_\times, \varepsilon_{\mathsf{pca}}, p \in (0,1)$. Then,* `LazySVD` *outputs a (column) orthonormal matrix $V_k = (v_1, \dots, v_k) \in \mathbb{R}^{d\times k}$ which, with probability at least $1-p$, satisfies all of the following properties.*

*(Denote by $M_k = (I - V_k V_k^\top) M (I - V_k V_k^\top)$.)*

*(a) Core lemma: if $\varepsilon_{\mathsf{pca}} \leq \frac{\varepsilon^4 \delta_\times^2}{2^{12} k^4 (\lambda_1/\lambda_k)^2}$, then $\|V_k^\top U\|_2 \leq \varepsilon$, where $U = (u_j, \dots, u_d)$ is the (column) orthonormal matrix and $j$ is the smallest index satisfying $\lambda_j \leq (1-\delta_\times)\|M_{k-1}\|_2$.*

*(b) Spectral norm guarantee: if $\varepsilon_{\mathsf{pca}} \leq \frac{\delta_\times^6}{2^{28} k^4 (\lambda_1/\lambda_{k+1})^6}$, then $\lambda_{k+1} \leq \|M_k\|_2 \leq \frac{\lambda_{k+1}}{1-\delta_\times}$.*

*(c) Rayleigh quotient guarantee: if $\varepsilon_{\mathsf{pca}} \leq \frac{\delta_\times^6}{2^{28} k^4 (\lambda_1/\lambda_{k+1})^6}$, then $(1-\delta_\times)\lambda_k \leq v_k^\top M v_k \leq \frac{1}{1-\delta_\times}\lambda_k$.*

*(d) Schatten-q norm guarantee: for every $q \geq 1$, if $\varepsilon_{\mathsf{pca}} \leq \frac{\delta_\times^6}{2^{28} k^4 d^{4/q} (\lambda_1/\lambda_{k+1})^6}$, then $\|M_k\|_{S_q} \leq \frac{(1+\delta_\times)^2}{(1-\delta_\times)^2} \Big( \sum_{i=k+1}^d \lambda_i^q \Big)^{1/q}$.*

*Proof of Theorem 4.1.* Let $V_s = (v_1, \dots, v_s)$, so we can write

$$M_s = (I - V_s V_s^\top) M (I - V_s V_s^\top) = (I - v_s v_s^\top) M_{s-1} (I - v_s v_s^\top)$$

(a) Define $\widehat{\lambda} = \|M_{k-1}\|_2 \geq \lambda_k$.

Note that all column vectors in $V_s$ are automatically eigenvectors of $M_s$ with eigenvalues zero. For analysis purpose only, let $W_s$ be the column matrix of eigenvectors in $V_s^\perp$ of $M_s$ that have eigenvalues in the range $[0, (1-\delta_\times+\tau_s)\widehat{\lambda}]$, where $\tau_s \stackrel{\text{def}}{=} \frac{s}{2k}\delta_\times$. We now show that for every $s = 0, \dots, k$, there exists a matrix $Q_s$ such that $\|U - W_s Q_s\|_2$ is small and $\|Q_s\|_2 \leq 1$. We will do this by induction.

In the base case: since $\tau_0 = 0$, we have $W_0 = U$ by the definition of $U$. We can therefore define $Q_0$ to be the identity matrix.

For every $s = 0, 1, \dots, k-1$, suppose there exists a matrix $Q_s$ with $\|Q_s\|_2 \leq 1$ that satisfies $\|U - W_s Q_s\|_2 \leq \eta_s$ for some $\eta_s > 0$, we construct $Q_{s+1}$ as follows.

First we observe that `AppxPCA` outputs a unit vector $v'_{s+1}$ satisfying $\|v'^\top_{s+1} W_s\|_2^2 \leq \varepsilon_{\mathsf{pca}}$ and $\|v'^\top_{s+1} V_s\|_2^2 \leq \varepsilon_{\mathsf{pca}}$ with probability at least $1 - p/k$. This follows from Theorem 3.1 because $[0, (1-\delta_\times+\tau_s)\widehat{\lambda}] \subseteq [0, (1-\delta_\times/2)\widehat{\lambda}]$, together with the fact that $\|M_s\|_2 \geq \|M_{k-1}\|_2 \geq \widehat{\lambda}$. Now, since $v_{s+1}$ is the projection of $v'_{s+1}$ into $V_s^\perp$, we have

$$\|v_{s+1}^\top W_s\|_2^2 \leq \frac{\|v'^\top_{s+1} W_s\|_2^2}{\|(I - V_s V_s^\top) v'_{s+1}\|_2^2} = \frac{\|v'^\top_{s+1} W_s\|_2^2}{1 - \|V_s^\top v'_{s+1}\|_2^2} \leq \frac{\varepsilon_{\mathsf{pca}}}{1 - \varepsilon_{\mathsf{pca}}} < 1.5\varepsilon_{\mathsf{pca}} \ . \tag{C.1}$$

Next we apply Lemma B.4 with $M = M_s$, $M' = M_{s+1}$, $U = W_s$, $V = W_{s+1}$, $v = v_{s+1}$, $\mu = (1-\delta_\times+\tau_s)\widehat{\lambda}$, and $\tau = (\tau_{s+1} - \tau_s)\widehat{\lambda}$. We obtain a matrix $\widetilde{Q}_s$, $\|\widetilde{Q}_s\|_2 \leq 1$ such that[7]

$$\|W_s - W_{s+1}\widetilde{Q_s}\|_2 \leq \sqrt{\frac{169(\lambda_1/\widehat{\lambda})^2 \cdot 1.5\varepsilon_{\mathsf{pca}}}{(\tau_{s+1}-\tau_s)^2} + \varepsilon_{\mathsf{pca}}} < \frac{32\lambda_1 k \sqrt{\varepsilon_{\mathsf{pca}}}}{\lambda_k \delta_\times} \ ,$$

[7]Technically speaking, to apply Lemma B.4 we need $U = W_s$ to consist of all eigenvectors of $M_s$ with eigenvalues $\leq \mu$. However, we only defined $W_s$ to be eigenvectors of $M_s$ with eigenvalues $\leq \mu$ that are *also* orthogonal to $V_s$. It is straightforward to verify that the same result of Lemma B.4 remains true because $v_{s+1}$ is orthogonal to $V_s$.

and this implies that

$$\|W_{s+1}\widetilde{Q_s}Q_s - U\|_2 \le \|W_{s+1}\widetilde{Q_s}Q_s - W_sQ_s\|_2 + \|W_sQ_s - U\|_2 \le \eta_s + \frac{32\lambda_1 k\sqrt{\varepsilon_{\mathsf{pca}}}}{\lambda_k\delta_\times} \enspace .$$

Let $Q_{s+1} = \widetilde{Q_s}Q_s$ we know that $\|Q_{s+1}\|_2 \le 1$ and

$$\|W_{s+1}Q_{s+1} - U\|_2 \le \eta_{s+1} \stackrel{\mathrm{def}}{=} \eta_s + \frac{32\lambda_1 k\sqrt{\varepsilon_{\mathsf{pca}}}}{\lambda_k\delta_\times} \enspace .$$

Therefore, after $k$-iterations of `LazySVD`, we obtain:

$$\|W_kQ_k - U\|_2 \le \eta_k = \frac{32\lambda_1 k^2\sqrt{\varepsilon_{\mathsf{pca}}}}{\lambda_k\delta_\times}$$

Multiply $U^\top$ from the left, we obtain $\|U^\top W_kQ_k - I\|_2 \le \eta_k$. Since $\|Q_k\|_2 \le 1$, we must have $\sigma_{\min}(U^\top W_k) \ge 1 - \eta_k$ (here $\sigma_{\min}$ denotes the smallest singular value). Therefore,

$$U^\top W_k W_k^\top U \succeq (1-\eta_k)^2 I \enspace .$$

Since $V_k$ and $W_k$ are orthogonal of each other, we have

$$U^\top V_k V_k^\top U \preceq U^\top (I - W_kW_k^\top)U \preceq I - (1-\eta_k)^2 I \preceq 2\eta_k I$$

Therefore,

$$\|V_k^\top U\|_2 \le 8\frac{(\lambda_1/\lambda_k)^{1/2} k \varepsilon_{\mathsf{pca}}^{1/4}}{\delta_\times^{1/2}} \le \varepsilon \enspace .$$

(b) The statement is obvious when $k = 0$. For every $k \ge 1$, the lower bound is obvious. We prove the upper bound by contradiction. Suppose that $\|M_k\|_2 > \frac{\lambda_{k+1}}{1-\delta_\times}$. Then, since $\|M_{k-1}\|_2 \ge \|M_k\|_2$ and therefore $\lambda_{k+1},\dots,\lambda_d < (1-\delta_\times)\|M_{k-1}\|_2$, we can apply Theorem 4.1.a of the current $k$ to deduce that $\|V_k^\top U_{>k}\|_2 \le \varepsilon$ where $U_{>k} \stackrel{\mathrm{def}}{=} (u_{k+1},\dots,u_d)$. We now apply Lemma B.2 with $V_s = V_k$ and $U = U_{>k}$, we obtain a matrix $Q_k \in \mathbb{R}^{d\times k}$ whose columns are spanned by $u_1,\dots,u_k$ and satisfy

$$\left\|\left(I - Q_kQ_k^\top\right) M \left(I - Q_kQ_k^\top\right) - \left(I - V_kV_k^\top\right) M \left(I - V_kV_k^\top\right)\right\|_2 < 16\lambda_1\varepsilon \enspace .$$

However, our assumption says that the second matrix $\left(I - V_kV_k^\top\right) M \left(I - V_kV_k^\top\right)$ has spectral norm at least $\lambda_{k+1}/(1-\delta_\times)$, but we know that $\left(I - Q_kQ_k^\top\right) M \left(I - Q_kQ_k^\top\right)$ has spectral norm exactly $\lambda_{k+1}$ due to the definition of $Q_k$. Therefore, we must have $\frac{\lambda_{k+1}}{1-\delta_\times} - \lambda_{k+1} \le 16\lambda_1\varepsilon$ due to triangle inequality.

In other words, by selecting $\varepsilon$ in Theorem 4.1.a to satisfy $\varepsilon \le \frac{\delta_\times}{16\lambda_1/\lambda_{k+1}}$ (which is satisfied by our assumption on $\varepsilon_{\mathsf{pca}}$), we get a contradiction so can conclude that $\|M_k\|_2 \le \frac{\lambda_{k+1}}{1-\delta_\times}$.

(c) We compute that

$$v_k^\top M v_k = v_k^\top M_{k-1} v_k \stackrel{①}{\ge} \frac{v_k'^\top M_{k-1} v_k'}{\|(I - V_{k-1}V_{k-1}^\top)v_k'\|_2^2} \stackrel{②}{\ge} \frac{v_k'^\top M_{k-1} v_k'}{1-\varepsilon_{\mathsf{pca}}}$$

$$\stackrel{③}{\ge} (1-\delta_\times/2)\|M_{k-1}\|_2 \ge (1-\delta_\times)\|M_{k-1}\|_2 \enspace .$$

Above, ① is because $v_k$ is the projection of $v'_k$ into $V_{k-1}^{\perp}$, ② is because $\|V_{k-1}^{\top} v'_k\|_2^2 \leq \varepsilon_{\mathsf{pca}}$ following the same reason as (C.1), and ③ is owing to Theorem 3.1. Next, since $\|M_{k-1}\|_2 \geq \lambda_k$, we automatically have $v_k^{\top} M v_k \geq (1-\delta_\times)\lambda_k$. On the other hand, $v_k^{\top} M v_k = v_k^{\top} M_{k-1} v_k \leq \|M_{k-1}\|_2 \leq \frac{\lambda_k}{1-\delta_\times}$ where the last inequality is owing to Theorem 4.1.b.

(d) Since $\|V_k^{\top} U\|_2 \leq \varepsilon_c \stackrel{\text{def}}{=} 8\frac{(\lambda_1/\lambda_k)^{1/2} k \varepsilon_{\mathsf{pca}}^{1/4}}{\delta_\times^{1/2}}$ from the analysis of Theorem 4.1.a, we can apply Lemma B.2 to obtain a (column) orthogonal matrix $Q_k \in \mathbb{R}^{d\times k}$ such that

$$\|M'_k - M_k\|_2 \leq 16\lambda_1\varepsilon_c, \qquad \text{where } M'_k \stackrel{\text{def}}{=} (I - Q_k Q_k^{\top}) M (I - Q_k Q_k^{\top}) \tag{C.2}$$

Suppose $U = (u_{d-p+1}, \dots, u_d)$ is of dimension $d \times p$, that is, there are exactly $p$ eigenvalues of $M$ that are $\leq (1-\delta_\times)\|M_{k-1}\|_2$. Then, the definition of $Q_k$ in Lemma B.2 tells us $U^{\top} Q_k = 0$ so $M'_k$ agrees with $M$ on all the eigenvalues and eigenvectors $\{(\lambda_j, u_j)\}_{j=d-p+1}^{d}$ because an index $j$ satisfies $\lambda_j \leq (1-\delta_\times)\|M_{k-1}\|_2$ if and only if $j \in \{d-p+1, d-p+2, \dots, d\}$.

Denote by $\mu_1, \dots, \mu_{d-k}$ the eigenvalues of $M'_k$ excluding the $k$ zero eigenvalues in subspace $Q_k$, and assume without loss of generality that $\{\mu_1, \dots, \mu_p\} = \{\lambda_{d-p+1}, \dots, \lambda_d\}$. Then,

$$\begin{aligned}
\|M'_k\|_{S_q}^q &= \sum_{i=1}^{d-k} \mu_i^q = \sum_{i=1}^{p} \mu_i^q + \sum_{i=p+1}^{d-k} \mu_i^q = \sum_{i=d-p+1}^{d} \lambda_i^q + \sum_{i=p+1}^{d-k} \mu_i^q \stackrel{①}{\leq} \sum_{i=d-p+1}^{d} \lambda_i^q + (d-k-p)\|M'_k\|_2^q \\
&\stackrel{②}{\leq} \sum_{i=d-p+1}^{d} \lambda_i^q + (d-k-p)(\|M_k\|_2 + 16\lambda_1\varepsilon_c)^q \\
&\stackrel{③}{\leq} \sum_{i=d-p+1}^{d} \lambda_i^q + (d-k-p)\Big(\frac{\lambda_{k+1}}{(1-\delta_\times)} + 16\lambda_1\varepsilon_c\Big)^q
\end{aligned}$$

Above, ① is because each $\mu_i$ is no greater than $\|M'_k\|_2$, and ② is owing to (C.2), and ③ is because of Theorem 4.1.b. Suppose we choose $\varepsilon_c$ so that $\varepsilon_c \leq \frac{\lambda_{k+1}\delta_\times}{16\lambda_1}$ (and this is indeed satisfied by our assumption on $\varepsilon_{\mathsf{pca}}$), then we can continue and write

$$\begin{aligned}
\|M'_k\|_{S_q}^q &\leq \sum_{i=d-p+1}^{d} \lambda_i^q + (d-k-p)\frac{(1+\delta_\times)^q}{(1-\delta_\times)^q}\lambda_{k+1}^q \\
&\stackrel{④}{\leq} \sum_{i=d-p+1}^{d} \lambda_i^q + \frac{(1+\delta_\times)^q}{(1-\delta_\times)^{2q}} \sum_{i=k+1}^{d-p} \lambda_i^q \leq \frac{(1+\delta_\times)^q}{(1-\delta_\times)^{2q}} \sum_{i=k+1}^{d} \lambda_i^q \ .
\end{aligned}$$

Above, ④ is because for each eigenvalue $\lambda_i$ where $i \in \{k+1, k+2, \dots, d-p\}$, we have $\lambda_i > (1-\delta_\times)\|M_{k-1}\|_2 \geq (1-\delta_\times)\lambda_k \geq (1-\delta_\times)\lambda_{k+1}$. Finally, using (C.2) again we have

$$\begin{aligned}
\|M_k\|_{S_q} &\leq \|M'_k\|_{S_q} + \|M_k - M'_k\|_{S_q} \leq \|M'_k\|_{S_q} + d^{1/p}\|M_k - M'_k\|_2 \\
&\leq \frac{1+\delta_\times}{(1-\delta_\times)^2}\Big(\sum_{i=k+1}^{d} \lambda_i^q\Big)^{1/q} + 16 d^{1/p}\lambda_1\varepsilon_c
\end{aligned}$$

As long as $\varepsilon_c \leq \frac{\delta_\times \lambda_{k+1}}{16 d^{1/p}\lambda_1}$, we have

$$\|M_k\|_{S_q} \leq \frac{(1+\delta_\times)^2}{(1-\delta_\times)^2}\Big(\sum_{i=k+1}^{d} \lambda_i^q\Big)^{1/q}$$

as desired. Finally, we note that $\varepsilon_c \leq \frac{\delta_\times \lambda_{k+1}}{16 d^{1/p}\lambda_1}$ is satisfied with our assumption on $\varepsilon_{\mathsf{pca}}$. □

### C.1 Proofs of Corollary 4.3 and Corollary 4.4

We first note that since `LazySVD` outputs $v_i$ one by one, although we have only stated Theorem 4.1 for the last iteration $k$, the claimed properties (a)-(d) hold for *all* intermediate iterations $s = 1, \dots, k$.

*Proof of Corollary 4.3 from Theorem 4.1.* By Theorem 4.1.a, we have: $\|V_k^\top U\|_2 \le \varepsilon$ where $U = (u_j, \dots, u_d)$ is a (column) orthonormal matrix and $j$ is the smallest index satisfying $\lambda_j \le (1 - \delta_\times)\|M_{k-1}\|_2$. Since it satisfies $\|M_{k-1}\|_2 \ge \lambda_k$, we have

$$\lambda_{k+1} = \sigma_{k+1}^2 = \sigma_k^2(1 - \mathsf{gap})^2 = \lambda_k(1 - \mathsf{gap})^2 \le \lambda_k(1 - \delta_\times) \le (1 - \delta_\times)\|M_{k-1}\|_2 \ ,$$

where the first inequality is because our choice of $\delta_\times = \mathsf{gap}$. Therefore, $j$ must be equal to $k+1$ according to its definition, and we conclude $\|V_k^\top W\|_2 \le \varepsilon$.

The running time of the algorithm comes directly from Theorem 4.2 by putting in the parameters. □

*Proof of Corollary 4.4 from Theorem 4.1.* Denote

$$M_k = (I - V_k V_k^\top) M (I - V_k V_k^\top) = (I - V_k V_k^\top) A A^\top (I - V_k V_k^\top) \ .$$

According to Theorem 4.1, we have:

$$\|M_k\|_{S_q} \le \left(\frac{1+\delta_\times}{1-\delta_\times}\right)^2 \|M - M_k^*\|_{S_q} \quad \forall q \ge 1 \text{ and } q = \infty \ ,$$

where $M_k^* = (A_k^*)(A_k^*)^\top$ is the rank-$k$ SVD of $M$, and recall $\delta_\times = \frac{\varepsilon}{3}$. For the spectral norm guarantee, we take $q = \infty$ and compute

$$\begin{aligned}
\|A - A_k\|_2 &= \|(I - V_k V_k^\top)A\|_2 = \sqrt{\|(I - V_k V_k^\top) A A^\top (I - V_k V_k^\top)\|_2} \\
&= \sqrt{\|M_k\|_2} = \sqrt{\|M_k\|_{S_\infty}} \le \frac{1+\delta_\times}{1-\delta_\times}\sqrt{\|M - M_k^*\|_{S_\infty}} \\
&\le (1+\varepsilon)\sqrt{\|M - M_k^*\|_{S_\infty}} = (1+\varepsilon)\sqrt{\|M - M_k^*\|_2} = (1+\varepsilon)\sigma_{k+1} = (1+\varepsilon)\|A - A_k^*\|_2 \ .
\end{aligned}$$

For the Frobenius norm guarantee, we take $q = 1$ and compute

$$\begin{aligned}
\|A - A_k\|_F &= \sqrt{\mathrm{Tr}[(I - V_k V_k^\top) A A^\top (I - V_k V_k^\top)]} = \sqrt{\|M_k\|_{S_1}} \le \frac{1+\delta_\times}{1-\delta_\times}\sqrt{\|M - M_k^*\|_{S_1}} \\
&\le (1+\varepsilon)\sqrt{\|M - M_k^*\|_{S_1}} = (1+\varepsilon)\sqrt{\sum_{i=k+1}^{d} \sigma_i^2} = (1+\varepsilon)\|A - A_k^*\|_F \ .
\end{aligned}$$

The Rayleigh quotient guarantees directly follow from Theorem 4.1.c. The running time of the algorithm comes directly from Theorem 4.2 by putting in the parameters. □

## D Proof Details for Our NNZ Running-Time Results

We state and prove a simple proposition first, and then divide this sections into three subsections: Section D.1 deals with column sampling together with the spectral-norm guarantee; Section D.2 deals with column sampling together with the Frobenius-norm guarantee; and Section D.3 deals with entry-wise sampling together with the spectral-norm guarantee.

**Proposition D.1.** *Let $A, A' \in \mathbb{R}^{d\times n}$ be two matrices with $d \le n$, and $\eta \ge 0$ be an non-negative real. Suppose $\|A - A'\|_2 \le \eta$, then for every $k \in [d]$, $\sigma_k(A') - \eta \le \sigma_k(A) \le \sigma_k(A') + \eta$*

*Proof of Proposition D.1.* By symmetry it is enough to show that $\sigma_k(A') \le \sigma_k(A) + \eta$.

Let $v_1, \dots, v_d$ be the (left) singular vectors of $A$ in decreasing order of the corresponding singular values, and let $S_k$ be the space spanned by $v_k, \dots, v_d$. Then, for every $x \in S_k$ that has $\|x\|_2 = 1$,

$$\|x^\top A'\|_2 \le \|x^\top (A - A')\|_2 + \|x^\top A\|_2 \le \sigma_k(A) + \eta \ .$$

Recall that the Courant-Fischer theorem says that

$$\sigma_k(A') = \min_{S, \dim(S) = d-k+1} \max_{x\in S, \|x\|_2 = 1} \|x^\top A'\|_2 \ .$$

Take $S = S_k$, we immediately obtain $\sigma_k(A') \le \sigma_k(A) + \eta$. □

## D.1 Column Sampling with Spectral-Norm Guarantee

We first state a concentration bound on column sampling (which is easily provable using for instance [22, Theorem 6.6.1]):

**Lemma D.2** (column sampling)**.** *Let $A \in \mathbb{R}^{d\times n}$ be a matrix and $A_i \in \mathbb{R}^d$ be the i-th column of A. Setting $p_i \stackrel{\mathrm{def}}{=} \|A_i\|_2^2 / \|A\|_F^2$ for each $i$, and define random rank-1 matrix $R = \frac{1}{p_i} A_i A_i^\top$ with probability $p_i$ for each $i$. For every $m \ge 1$, define $\overline{R}_m$ to be the average of $m$ independent copies of $R$, that is, $\overline{R}_m \stackrel{\mathrm{def}}{=} \frac{1}{m}\sum_{t=1}^m R_t$ where each $R_t$ is drawn from $R$. Then, for every $\eta, \delta > 0$, if*

$$m \ge \frac{8\|A\|_F^2 \|A\|_2^2 \log\frac{1}{\delta}}{\eta^2} + \frac{8\|A\|_F^2 \log\frac{1}{\delta}}{\eta} \ ,$$

*we have that with probability $1 - \delta$, it satisfies $\left\|\overline{R}_m - AA^\top\right\|_2 \le \eta$.*

The next lemma translates the approximate solution on the column sampled matrix into a spectral guarantee on the original matrix.

**Lemma D.3.** *Let $A \in \mathbb{R}^{d\times n}$ be a matrix and define $\overline{R}_m$ as in Lemma D.2. For every $k \in [d-1]$, every $\varepsilon > 0$, every $p \in (0,1)$, and every $\delta_\times \in (0,1)$, if $m \ge \frac{32k\log(1/p)\sigma_1(A)^4}{\varepsilon^2 \sigma_{k+1}(A)^4}$ and one obtains an $\delta_\times$-approximate k-SVD of $\overline{R}_m$ in terms of spectral norm, that is*

$$\text{a column orthogonal matrix } V \in \mathbb{R}^{d\times k} \text{ such that } \|(I - VV^\top)\overline{R}_m(I - VV^\top)\|_2 \le \frac{\lambda_{k+1}(\overline{R}_m)}{1-\delta_\times} \ .$$

*Then, with probability at least $1 - p$, this matrix $V$ satisfies*

$$\|A - VV^\top A\|_2 \le \frac{\|A - A_k^*\|_2 + 2\varepsilon\|A - A_k^*\|_F}{1 - \delta_\times} \ .$$

*Proof of Lemma D.3.* If we let $\eta \stackrel{\text{def}}{=} \varepsilon^2\|A-A_k^*\|_F^2 + \varepsilon\|A-A_k^*\|_F\|A-A_k^*\|_2$, we can compute

$$\begin{aligned}
& \frac{8\|A\|_F^2\|A\|_2^2\log\frac{1}{p}}{\eta^2} + \frac{8\|A\|_F^2\log\frac{1}{p}}{\eta} \\
\leq\ & \frac{8\|A\|_F^2\|A\|_2^2\log\frac{1}{p}}{\varepsilon^2\|A-A_k^*\|_F^2\|A-A_k^*\|_2^2} + \frac{8\|A\|_F^2\log\frac{1}{p}}{\varepsilon^2\|A-A_k^*\|_F^2} \\
\leq\ & \frac{8\left(\|A-A_k^*\|_F^2 + k\|A\|_2^2\right)\|A\|_2^2\log\frac{1}{p}}{\varepsilon^2\|A-A_k^*\|_F^2\|A-A_k^*\|_2^2} + \frac{8\left(\|A-A_k^*\|_F^2 + k\|A\|_2^2\right)\log\frac{1}{p}}{\varepsilon^2\|A-A_k^*\|_F^2} \\
\leq\ & \frac{8\|A\|_2^2\log\frac{1}{p}}{\varepsilon^2\|A-A_k^*\|_2^2} + \frac{8k\|A\|_2^4\log\frac{1}{p}}{\varepsilon^2\|A-A_k^*\|_2^4} + \frac{8\log\frac{1}{p}}{\varepsilon^2} + \frac{8k\|A\|_2^2\log\frac{1}{p}}{\varepsilon^2\|A-A_k^*\|_2^2} \\
\leq\ & \frac{32k\log(1/p)\sigma_1(A)^4}{\varepsilon^2\sigma_{k+1}(A)^4} \leq m
\end{aligned}$$

Therefore, according to Lemma D.2, with probability at least $1-p$, it satisfies

$$\|AA^\top - \overline{R}_m\|_2 \leq \eta = \varepsilon^2\|A-A_k^*\|_F^2 + \varepsilon\|A-A_k^*\|_F\|A-A_k^*\|_2 \ .$$

This further implies that, owing to $\|I - VV^\top\|_2 \leq 1$,

$$\begin{aligned}
\|(I-VV^\top)AA^\top(I-VV^\top)\|_2 &\leq \|(I-VV^\top)\overline{R}_m(I-VV^\top)\|_2 + \|AA^\top - \overline{R}_m\|_2 \\
&\leq \frac{\lambda_{k+1}(\overline{R}_m)}{1-\delta_\times} + \eta \leq \frac{\lambda_{k+1}(AA^\top) + 2\eta}{1-\delta_\times} \\
&= \frac{\sigma_{k+1}(A)^2 + 2\varepsilon^2\|A-A_k^*\|_F^2 + 2\varepsilon\|A-A_k^*\|_F\|A-A_k^*\|_2}{1-\delta_\times} \\
&\leq \frac{(\|A-A_k^*\|_2 + 2\varepsilon\|A-A_k^*\|_F)^2}{1-\delta_\times} \ .
\end{aligned}$$

Therefore,

$$\|(I-VV^\top)A\|_2 = \sqrt{\|(I-VV^\top)AA^\top(I-VV^\top)\|_2} \leq \frac{\|A-A_k^*\|_2 + 2\varepsilon\|A-A_k^*\|_F}{1-\delta_\times} \ . \qquad \square$$

Using the previous lemma, it is not hard to deduce our main result of this sub-section.

**Theorem D.4.** *Let $A \in \mathbb{R}^{d\times n}$ be a matrix with singular values $\sigma_1 \geq \cdots \geq \sigma_d \geq 0$. For every $\varepsilon \in (0, 1/2)$, let $\overline{R}_m$ be the subsampled version of $A$ as defined in Lemma D.2 with $m = \Omega\big(\frac{k\log(1/\varepsilon)\sigma_1^4}{\varepsilon^2\sigma_{k+1}^4}\big)$. Then, one can call* `LazySVD` *with appropriately chosen $\delta_\times$ to produce a matrix $V_k \in \mathbb{R}^{d\times k}$ satisfying*

$$\|A - V_kV_k^\top A\|_2 \leq \|A-A_k^*\|_2 + O(\varepsilon)\|A-A_k^*\|_F \ , \tag{D.1}$$

*and the total running time is $O(\mathsf{nnz}(A)) + \widetilde{O}\big(\frac{k^2 d(\sigma_1/\sigma_{k+1})^4}{\varepsilon^{2.5}}\big)$ if AGD is used as the approximate matrix inversion algorithm $\mathcal{A}$. Furthermore, if $\varepsilon \leq \frac{\|A-A_k^*\|_2}{2\|A-A_k^*\|_F}$, then one can use accelerated SVRG as $\mathcal{A}$ and improve the running time to $O(\mathsf{nnz}(A)) + \widetilde{O}\big(\frac{k^2 d(\sigma_1/\sigma_{k+1})^4}{\varepsilon^2}\big)$.*

*Proof of Theorem D.4.* Let us define $\overline{M} \stackrel{\text{def}}{=} \overline{R}_m/\|A\|_F^2$, and we can write $\overline{M} = \frac{1}{m}\sum_{t=1}^m a_i a_i^\top$ where each $a_i$ has Euclidean norm at most 1 due to our definition of $\overline{R}_m$ in Lemma D.2. We pass this matrix $\overline{M}$ as the input $M$ to `LazySVD`, and it satisfies $\|\overline{M}\|_2 \leq \mathrm{Tr}(\overline{M}) \leq 1$.

Before specifying the parameter choices for $\delta_\times$ and the algorithm choice for $\mathcal{A}$, we notice that for sufficiently small $\varepsilon_{\mathsf{pca}}$, Theorem 4.1.b implies $\|(I - V_k V_k^\top)\overline{R}_m(I - V_k V_k^\top)\|_2 \leq \frac{1}{1-\delta_\times}\lambda_{k+1}(\overline{R}_m)$ where $V_k$ is the output matrix from `LazySVD`. Applying Lemma D.3, we have

$$\|A - V_k V_k^\top A\|_2 \leq \frac{\|A - A_k^*\|_2 + 2\varepsilon\|A - A_k^*\|_F}{1 - \delta_\times} \ .$$

Now there are two cases. If we use AGD as the method $\mathcal{A}$, then we can choose $\delta_\times = \varepsilon$. In such a case it is easy to see that $\frac{\|A-A_k^*\|_2+2\varepsilon\|A-A_k^*\|_F}{1-\delta_\times} \leq \|A - A_k^*\|_2 + 4\varepsilon\|A - A_k^*\|_F$ so the guarantee (D.1) is satisfied. The total running time is $O(\mathsf{nnz}(A))$ to sample $\{a_1, \dots, a_m\}$ plus $\widetilde{O}(kmd/\sqrt{\varepsilon})$ to perform `LazySVD` (see Theorem 4.2).

If we use accelerated SVRG as the method $\mathcal{A}$, then we obtain a better dependence on $\varepsilon$ as follows. Suppose $\varepsilon \leq \frac{\|A-A_k^*\|_2}{2\|A-A_k^*\|_F}$, and we choose $\delta_\times \leq \varepsilon\frac{\|A-A_k^*\|_F}{\|A-A_k^*\|_2} \leq 1/2$. In such a case, one can verify again that $\frac{\|A-A_k^*\|_2+2\varepsilon\|A-A_k^*\|_F}{1-\delta_\times} \leq \|A - A_k^*\|_2 + O(\varepsilon)\|A - A_k^*\|_F$ so the guarantee (D.1) is satisfied. As for the running time, in addition to $O(\mathsf{nnz}(A))$ for sampling, we need (using Theorem 4.2 again)

$$\widetilde{O}\Big(kmd + \frac{km^{3/4}d}{\lambda_k(\overline{M})^{1/4}\delta_\times^{1/2}}\Big) = \widetilde{O}\Big(kmd + \frac{km^{3/4}d}{\lambda_k(\overline{M})^{1/4}\delta_\times^{1/2}}\Big) = \widetilde{O}\Big(kmd + \frac{km^{3/4}d}{\delta_\times^{1/2}}\Big(\frac{\|A\|_F^2}{\sigma_k(A)^2}\Big)^{1/4}\Big)$$
$$= \widetilde{O}\Big(kmd + \frac{km^{3/4}d}{\varepsilon^{1/2}}\Big(\frac{\sigma_{k+1}(A)^2\|A\|_F^2}{\sigma_k(A)^2\|A - A_k^*\|_F^2}\Big)^{1/4}\Big)$$

Since it can be verified (similar to the proof of Lemma D.3) $\frac{\sigma_{k+1}(A)^2\|A\|_F^2}{\sigma_k(A)^2\|A-A_k^*\|_F^2} \leq 1 + \frac{\sigma_1(A)^2}{\sigma_k(A)^2} \cdot k$, we conclude that the above running time becomes

$$\widetilde{O}\Big(kmd + \frac{km^{3/4}d}{\varepsilon^{1/2}}\Big(\frac{\sigma_1(A)^2}{\sigma_k(A)^2} \cdot k\Big)^{1/4}\Big) = \widetilde{O}\Big(\frac{k^2 d(\sigma_1/\sigma_{k+1})^4}{\varepsilon^2}\Big)$$

where the last equality follows from the definition of $m$. □

## D.2 Column Sampling with Frobenius-Norm Guarantee

The next lemma translates the approximate solution on the column sampled matrix into a Frobenius-norm guarantee on the original matrix.

**Lemma D.5.** *Let $A \in \mathbb{R}^{d\times n}$ be a matrix and define $\overline{R}_m$ as in Lemma D.2. For every $k \in [d-1]$, every $\varepsilon > 0$, every $p \in (0,1)$, and every $\delta_\times \in (0,1)$, if $m \geq \frac{128k^3\log(1/p)\sigma_1(A)^4}{\varepsilon^2\sigma_{k+1}(A)^4}$ and one obtains a $\delta_\times$-approximate $k$-SVD of $\overline{R}_m$ in terms of Rayleigh quotient, that is, a column orthonormal matrix $V = (v_1, \dots, v_k) \in \mathbb{R}^{d\times k}$ such that*

$$\forall i \in [k], \quad |v_i^\top \overline{R}_m v_i - \lambda_i(\overline{R}_m)| \leq \delta_\times \lambda_i(\overline{R}_m) \ .$$

*Then, with probability at least $1 - p$, this matrix $V$ satisfies*

$$\|(I - VV^\top)A\|_F^2 \leq \delta_\times\|A_k^*\|_F^2 + (1 + \varepsilon)\|A - A_k^*\|_F^2$$

*Proof of Lemma D.5.* Using similar arguments as in the proof of Lemma D.3, one can deduce that our choice of $m$ ensures

$$\|AA^\top - \overline{R}_m\|_2 \leq \eta \stackrel{\text{def}}{=} \frac{\varepsilon}{2k}\|A - A_k^*\|_F^2 \ .$$

Which implies that $|\lambda_i(\overline{R}_m) - \sigma_i^2| \le \eta$ due to Proposition D.1. Finally,

$$\|(I - VV^\top)A\|_F^2 = \mathrm{Tr}[(I - VV^\top)AA^\top(I - VV^\top)] = \mathrm{Tr}[AA^\top - VAA^\top V]$$
$$= \sum_{i=1}^{d} \sigma_i(A)^2 - \sum_{i=1}^{k} v_i^\top AA^\top v_i \le \sum_{i=1}^{d} \sigma_i(A)^2 + k\eta - \sum_{i=1}^{k} v_i^\top \overline{R}_m v_i$$
$$\le \sum_{i=1}^{d} \sigma_i(A)^2 + k\eta - \sum_{i=1}^{k}(1-\delta_\times)(\sigma_i^2(A) - \eta) \le \sum_{i=1}^{k} \delta_\times \sigma_i(A)^2 + 2k\eta + \sum_{i=k+1}^{d} \sigma_i(A)^2$$
$$= \delta_\times \|A_k^*\|_F^2 + \|A - A_k^*\|_F^2 + 2k\eta = \delta_\times \|A_k^*\|_F^2 + (1+\varepsilon)\|A - A_k^*\|_F^2 \ . \qquad \square$$

Using the previous lemma, it is not hard to deduce our main result of this sub-section.

**Theorem D.6.** *Let $A \in \mathbb{R}^{d\times n}$ be a matrix with singular values $\sigma_1 \ge \dots \ge \sigma_d \ge 0$. For every $\varepsilon \in (0, 1/2)$, let $\overline{R}_m$ be the subsampled version of $A$ as defined in Lemma D.2 with $m = \Omega\big(\frac{k^3 \log(1/\varepsilon)\sigma_1^4}{\varepsilon^2 \sigma_{k+1}^4}\big)$. Then, one can call* `LazySVD` *with appropriately chosen $\delta_\times$ to produce a matrix $V_k \in \mathbb{R}^{d\times k}$ satisfying*

$$\|A - V_k V_k^\top A\|_2 \le (1 + O(\varepsilon))\|A - A_k^*\|_F \ ,$$

*and the total running time is $O(\mathsf{nnz}(A)) + \widetilde{O}\big(\frac{k^4 d(\sigma_1/\sigma_{k+1})^5}{\varepsilon^{2.5}}\big)$ if AGD is used as the approximate matrix inversion algorithm $\mathcal{A}$, or $O(\mathsf{nnz}(A)) + \widetilde{O}\big(\frac{k^4 d(\sigma_1/\sigma_{k+1})^{4.5}}{\varepsilon^2}\big)$ if $\mathcal{A}$ is accelerated SVRG.*

*Proof.* One can choose $\delta_\times = \varepsilon\|A - A_k^*\|_F^2/\|A_k^*\|_F^2 \ge \frac{\varepsilon}{k}\frac{\sigma_{k+1}^2}{\sigma_1^2}$ as the parameter of `LazySVD` and the proof is completely analogous to that of Theorem D.4. $\square$

## D.3 Entry-Wise Sampling with Spectral-Norm Guarantee

We first state a concentration bound on entry-wise sampling:

**Lemma D.7** (entry-wise sampling [10])**.** *Let $A \in \mathbb{R}^{d\times n}$ be a matrix. Define random single-entry matrix $R = \frac{1}{p_{i,j}} A_{i,j}$ where $(i,j)$ is selected from $[d] \times [m]$ each with probability $p_{i,j} \overset{\text{def}}{=} A_{i,j}^2/\|A\|_F^2$. For every $m \ge 1$, define $\overline{R}_s$ to be the average of $s$ independent copies of $R$, that is, $\overline{R}_s = \frac{1}{s}\sum_{t=1}^{s} R_t$ where each $R_t$ is drawn from $R$. Then, for every $\eta, p \in (0,1)$, if*

$$s \ge \frac{28(d+n)\log(2/p)\|A\|_F^2}{\eta^2} \ ,$$

*we have with probability $1-p$, it satisfies $\|A - \overline{R}_s\|_2 \le \eta$.*

The next lemma translates the approximate solution on the entry-wise sampled matrix into a spectral guarantee on the original matrix.

**Lemma D.8.** *Let $A \in \mathbb{R}^{d\times n}$ be a matrix and define $\overline{R}_s$ as in Lemma D.7. For every $k \in [d-1]$, every $\varepsilon > 0$, every $p \in (0,1]$, and every $\delta_\times \in (0,1)$, if $s \ge \frac{56k(d+n)\ln(2/p)\sigma_1(A)^2}{\varepsilon^2 \sigma_{k+1}(A)^2}$ and one obtains a matrix $R'$ satisfying*

$$\|R' - \overline{R}_s\|_2 \le \frac{\sigma_{k+1}(\overline{R}_s)}{1-\delta_\times} \ .$$

*Then, with probability at least $1-p$, we also have*

$$\|A - R'\|_2 \le \frac{\|A - A_k^*\|_2 + 2\varepsilon\|A - A_k^*\|_F}{1-\delta_\times} \ .$$

*Proof of Lemma D.8.* We first compute that

$$
\begin{aligned}
s &\geq \frac{56k(d+n)\ln(2/p)\sigma_1(A)^2}{\varepsilon^2\sigma_{k+1}(A)^2} \\
&\geq \frac{28k(d+n)\ln(2/p)\sigma_1(A)^2}{\varepsilon^2\sigma_{k+1}(A)^2} + \frac{28d\ln(2/p)}{\varepsilon^2} \\
&\geq \frac{28(d+n)\ln(2/p)\|A_k^*\|_F^2}{\varepsilon^2\|A-A_k^*\|_F^2} + \frac{28(d+n)\ln(2/p)\|A-A_k^*\|_F^2}{\varepsilon^2\|A-A_k^*\|_F^2} \\
&= \frac{28(d+n)\ln(2/p)\left(\|A_k^*\|_F^2 + \|A-A_k^*\|_F^2\right)}{\varepsilon^2\|A-A_k^*\|_F^2} \\
&= \frac{28(d+n)\ln(2/p)\|A\|_F^2}{\varepsilon^2\|A-A_k^*\|_F^2} \enspace.
\end{aligned}
$$

Owing to Lemma D.7, with probability $1-p$, it satisfies $\|A-\overline{R}_s\|_2 \leq \eta \stackrel{\text{def}}{=} \varepsilon\|A-A_k^*\|_F$. Using Proposition D.1, we know that $\sigma_{k+1}(\overline{R}_s) \leq \sigma_{k+1}(A)+\eta$, which implies that

$$\|R'-\overline{R}_s\|_2 \leq \frac{\sigma_{k+1}(A)+\eta}{1-\delta_\times} \enspace.$$

Finally,

$$\|R'-A\|_2 \leq \|R'-\overline{R}_s\|_2 + \|A-\overline{R}_s\|_2 \leq \frac{\sigma_{k+1}(A)+\eta}{1-\delta_\times} + \eta \leq \frac{\sigma_{k+1}(A)+2\varepsilon\|A-A_k^*\|_F}{1-\delta_\times} \enspace. \qquad \square$$

Using the previous lemma, it is not hard to deduce our main result of this sub-section.

**Theorem D.9.** *Let $A \in \mathbb{R}^{d\times n}$ be a matrix with singular values $\sigma_1 \geq \cdots \geq \sigma_d \geq 0$. For every $\varepsilon \in (0,1/2)$, let $\overline{R}_s$ be the entry-sampled version of $A$ as defined in Lemma D.7 with $s = \Omega\big(\frac{k(d+n)\log(1/\varepsilon)\sigma_1^2}{\varepsilon^2\sigma_{k+1}^2}\big)$. Then, one can call* `LazySVD` *with appropriately chosen $\delta_\times$ to produce a matrix $V_k \in \mathbb{R}^{d\times k}$ satisfying*

$$\|A - V_kV_k^\top A\|_2 \leq \|A-A_k^*\|_2 + O(\varepsilon)\|A-A_k^*\|_F \enspace,$$

*and the total running time is $O(\mathsf{nnz}(A)) + \widetilde{O}\big(\frac{k^2(n+d)(\sigma_1/\sigma_{k+1})^2}{\varepsilon^{2.5}}\big)$ if AGD is used as the approximate matrix inversion algorithm $\mathcal{A}$.*

*Proof.* One can choose $\delta_\times = \varepsilon$ and a completely analogous proof as that of Theorem D.4 gives us the desired spectral guarantee. The running time follows from Theorem 4.2 because $\mathsf{nnz}(\overline{R}_s) = s$. $\square$

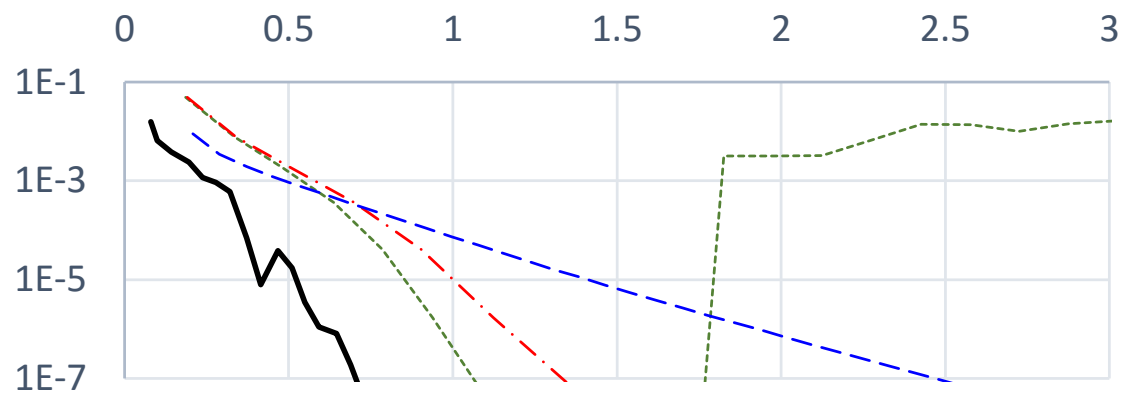

(a) email, $k = 10$, `Fnorm` metric

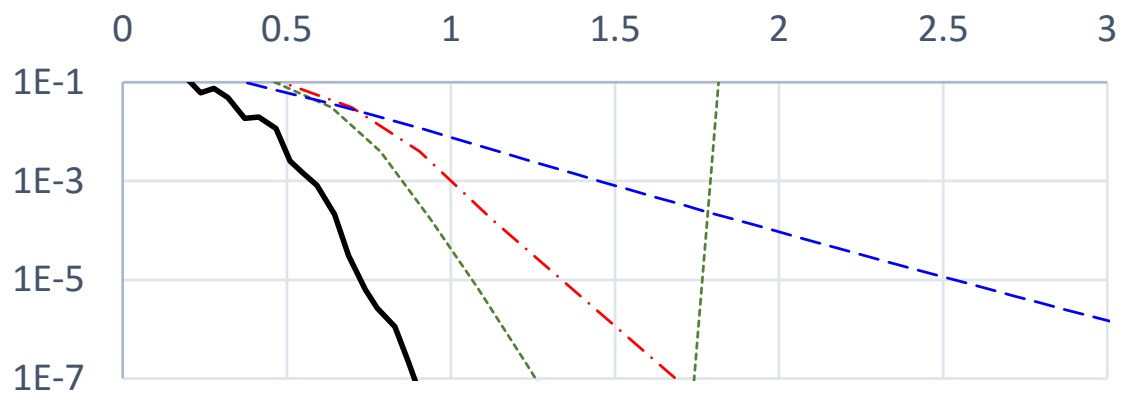

(b) email, $k = 10$, `rayleigh` metric

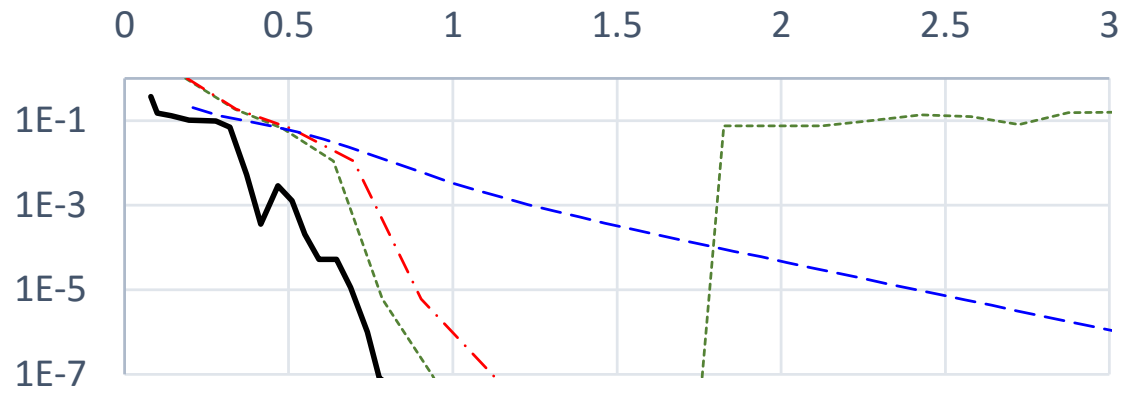

(c) email, $k = 10$, `spectral` metric

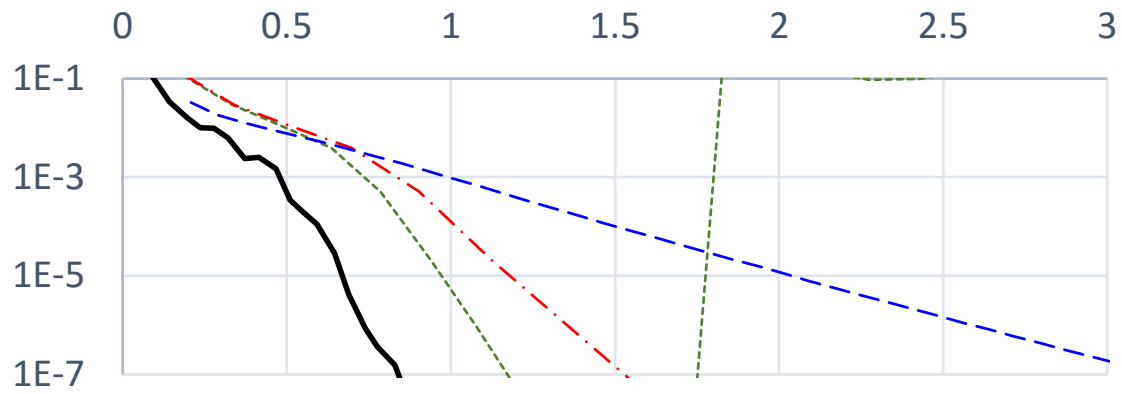

(d) email, $k = 10$, `rayleigh(last)` metric

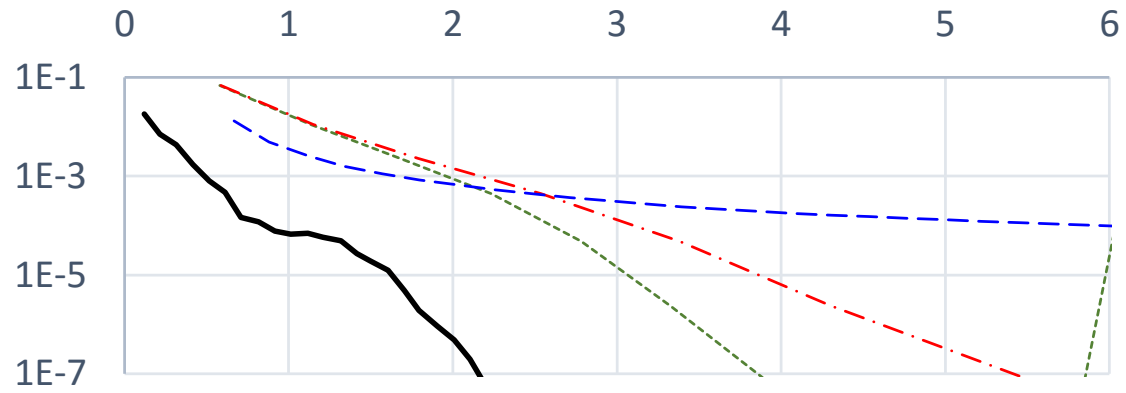

(e) email, $k = 20$, `Fnorm` metric

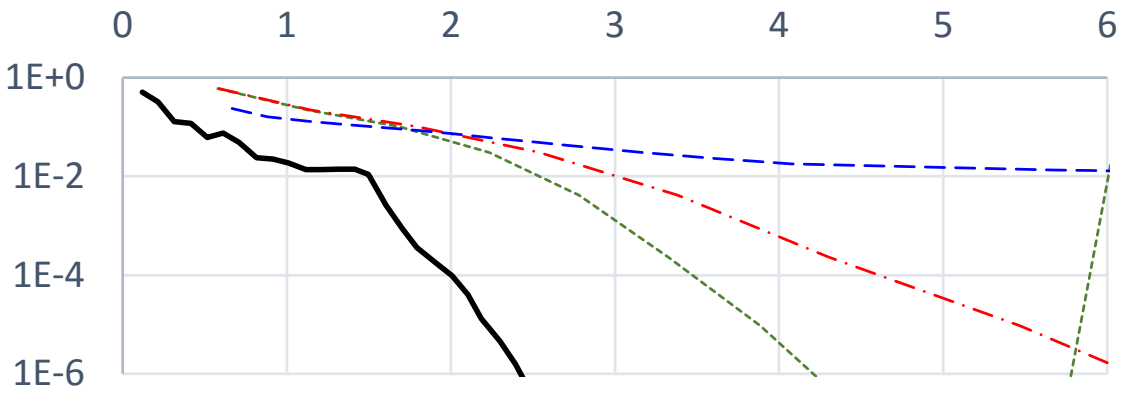

(f) email, $k = 20$, `rayleigh` metric

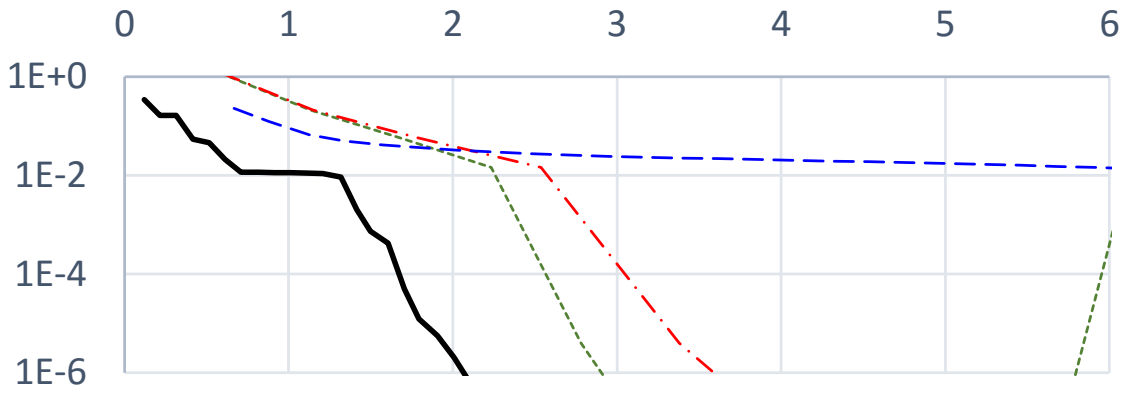

(g) email, $k = 20$, `spectral` metric

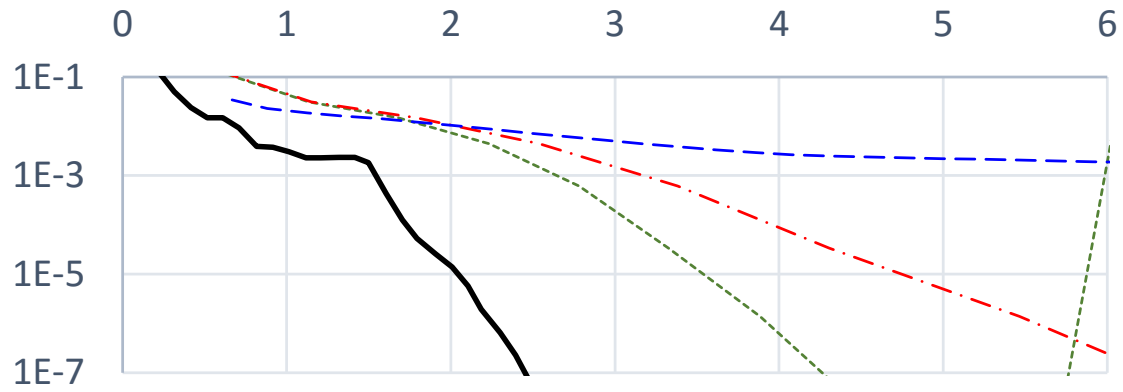

(h) email, $k = 20$, `rayleigh(last)` metric

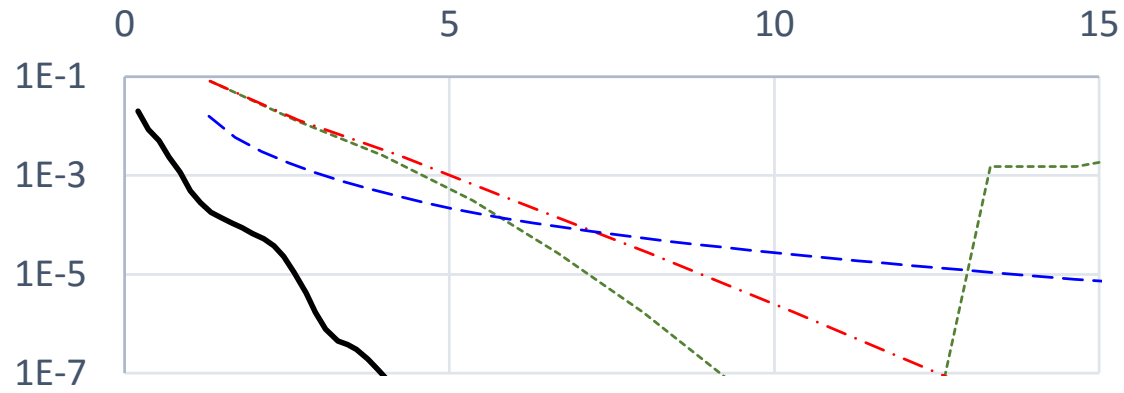

(i) email, $k = 30$, `Fnorm` metric

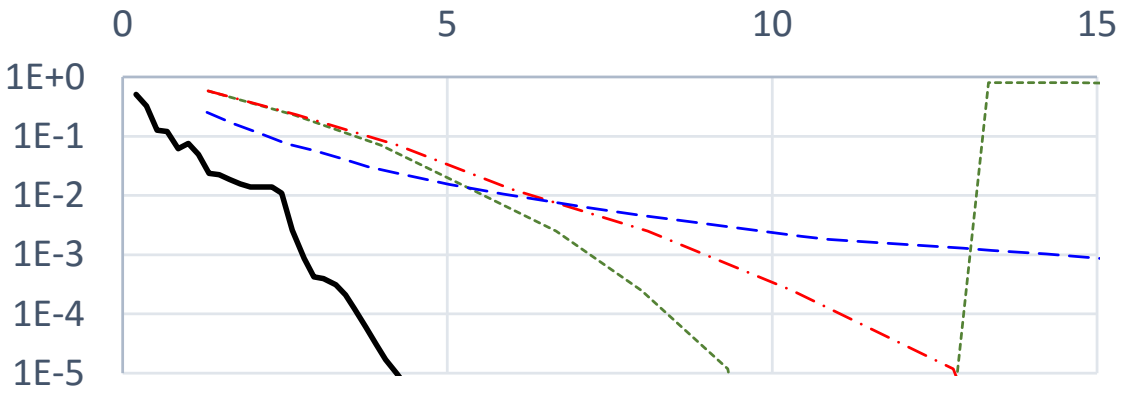

(j) email, $k = 30$, `rayleigh` metric

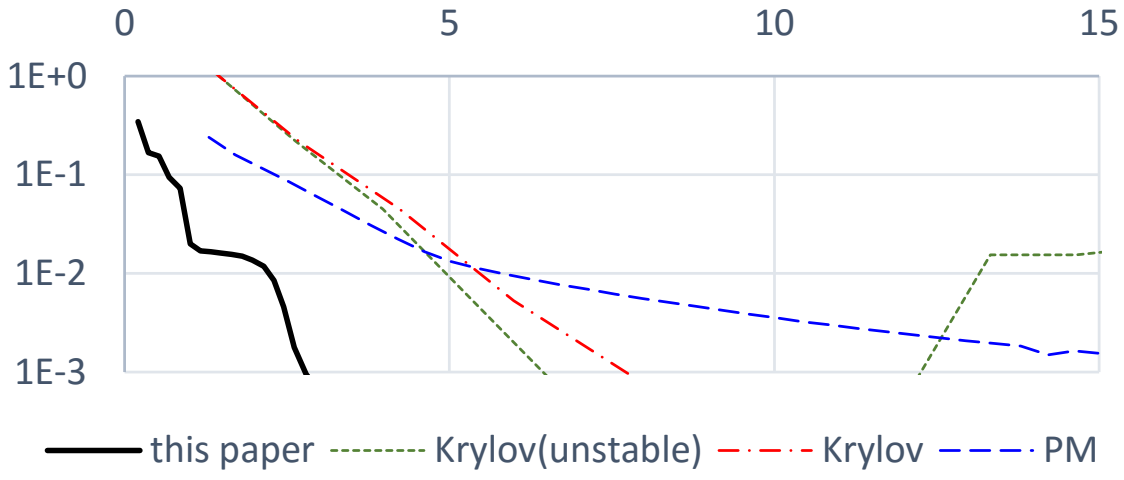

(k) email, $k = 30$, `spectral` metric

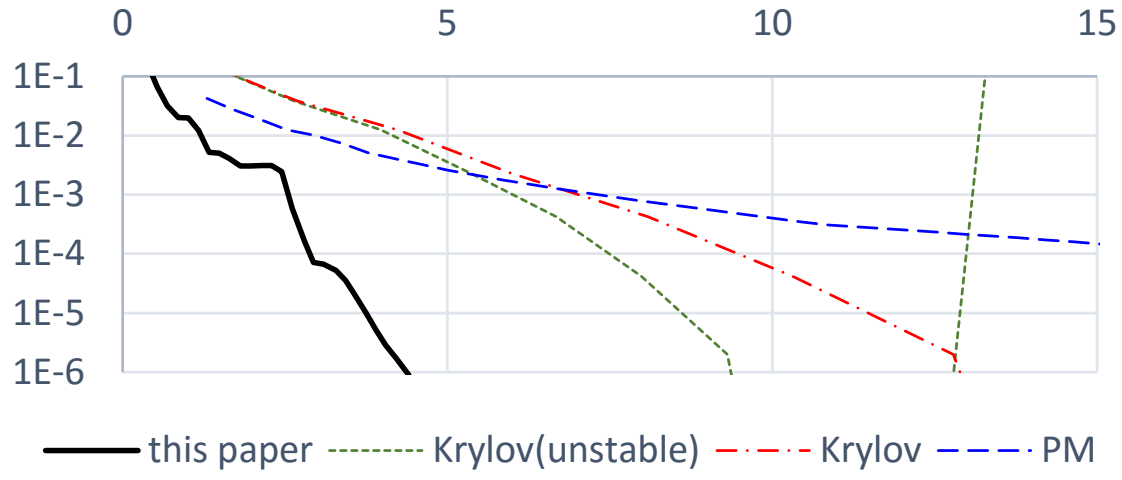

(l) email, $k = 30$, `rayleigh(last)` metric

Figure 2: Performance on dataset `email`. Relative error ($y$-axis) vs. running time ($x$-axis).

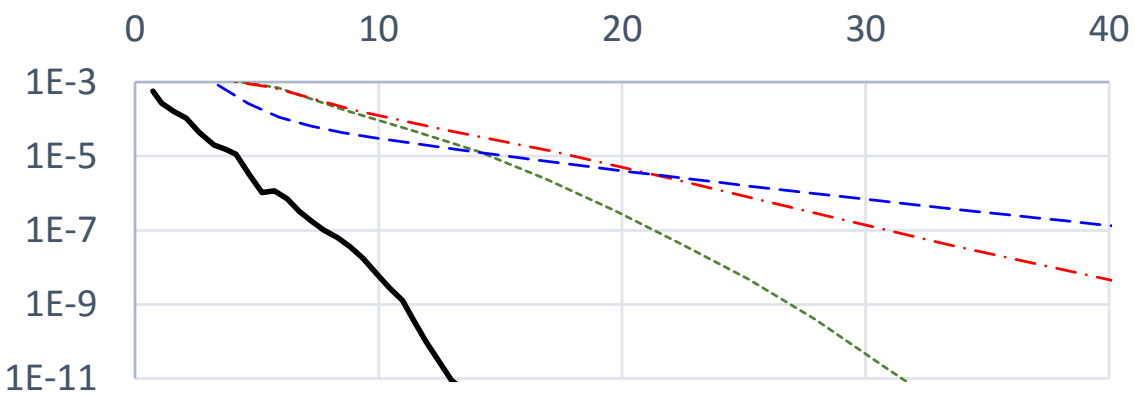

(a) amazon, $k = 10$, `Fnorm` metric

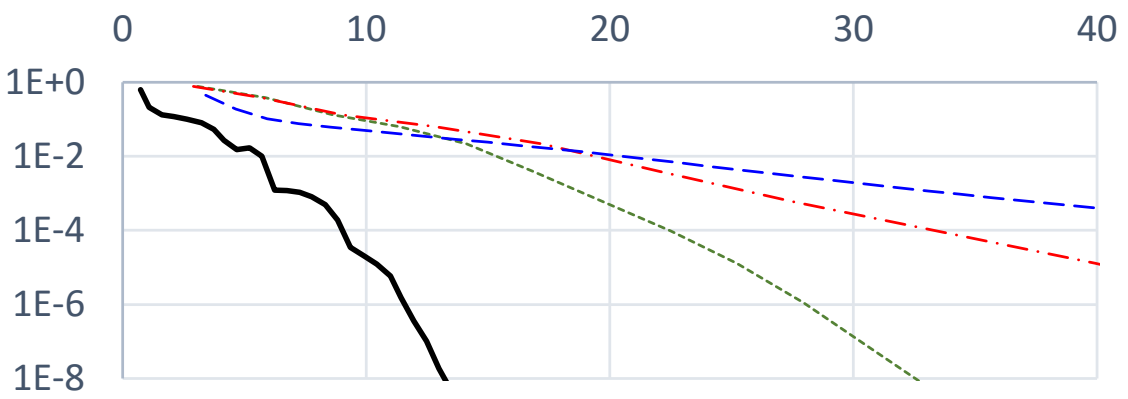

(b) amazon, $k = 10$, `rayleigh` metric

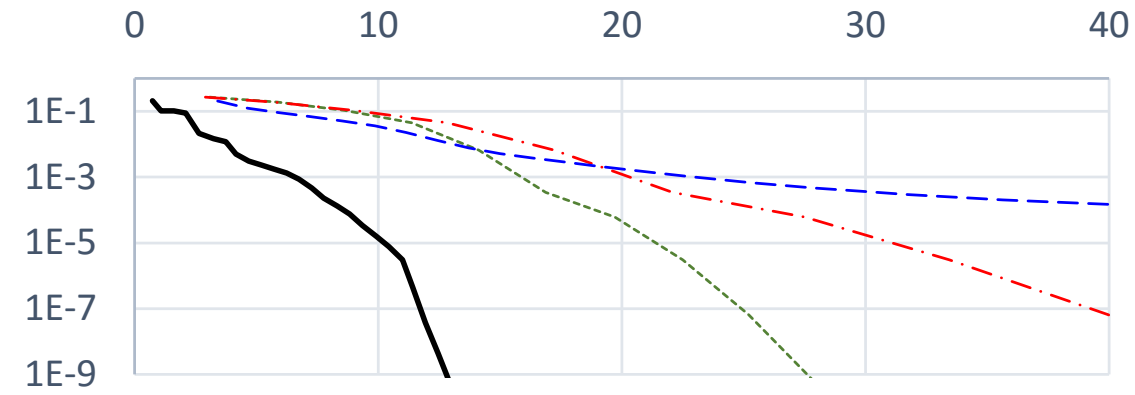

(c) amazon, $k = 10$, `spectral` metric

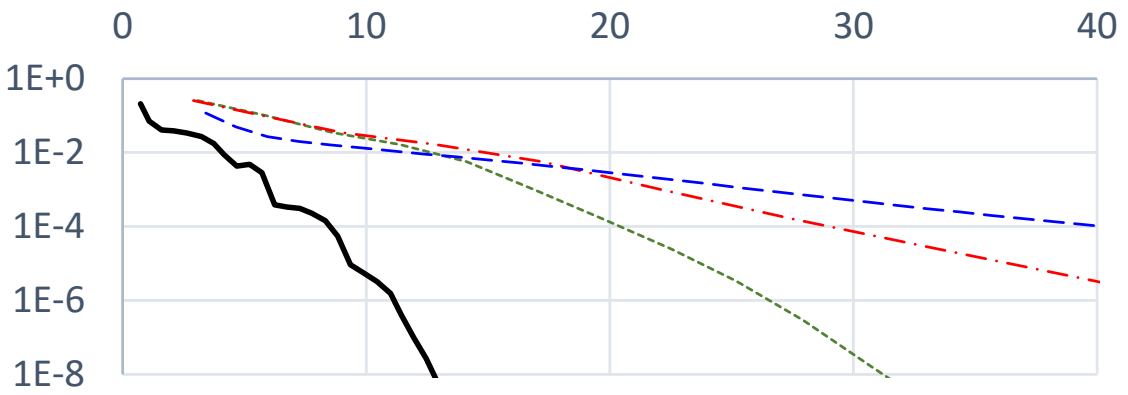

(d) amazon, $k = 10$, `rayleigh(last)` metric

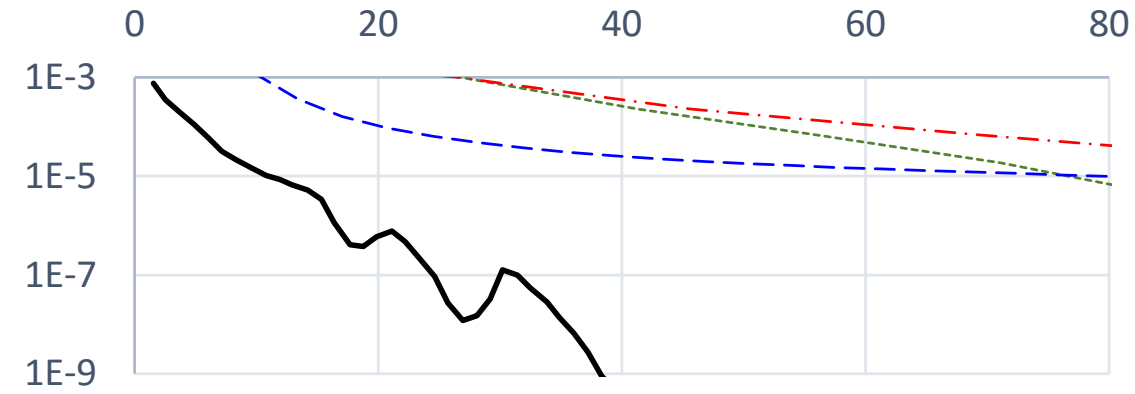

(e) amazon, $k = 20$, `Fnorm` metric

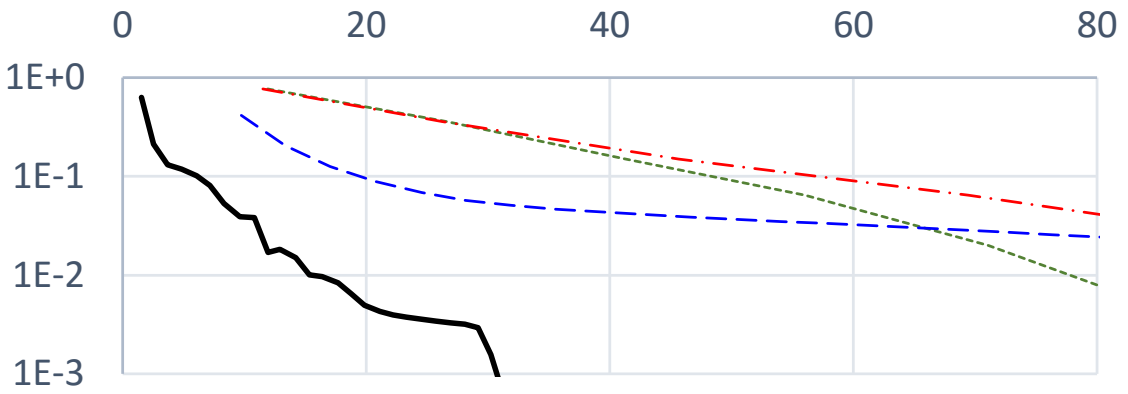

(f) amazon, $k = 20$, `rayleigh` metric

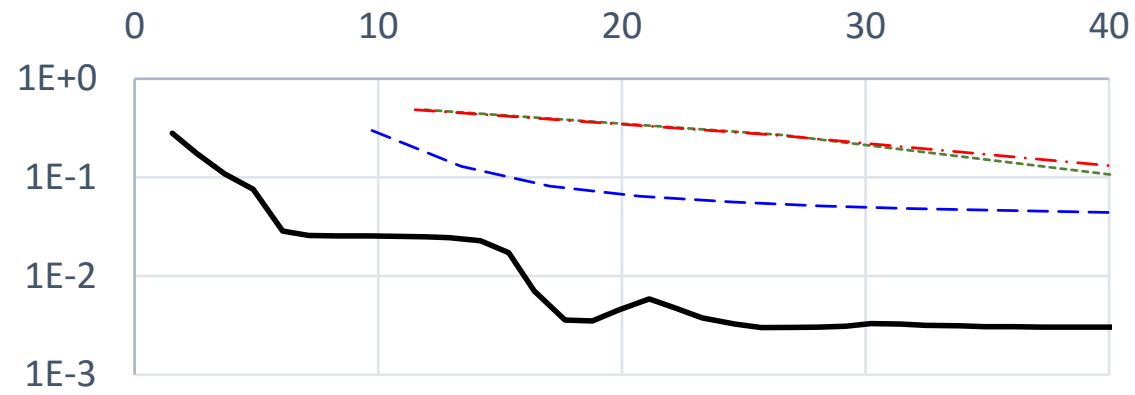

(g) amazon, $k = 20$, `spectral` metric

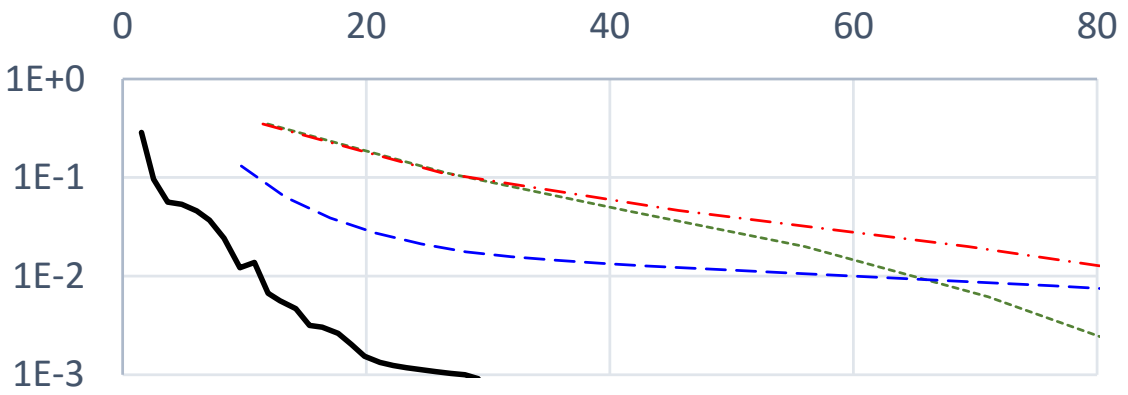

(h) amazon, $k = 20$, `rayleigh(last)` metric

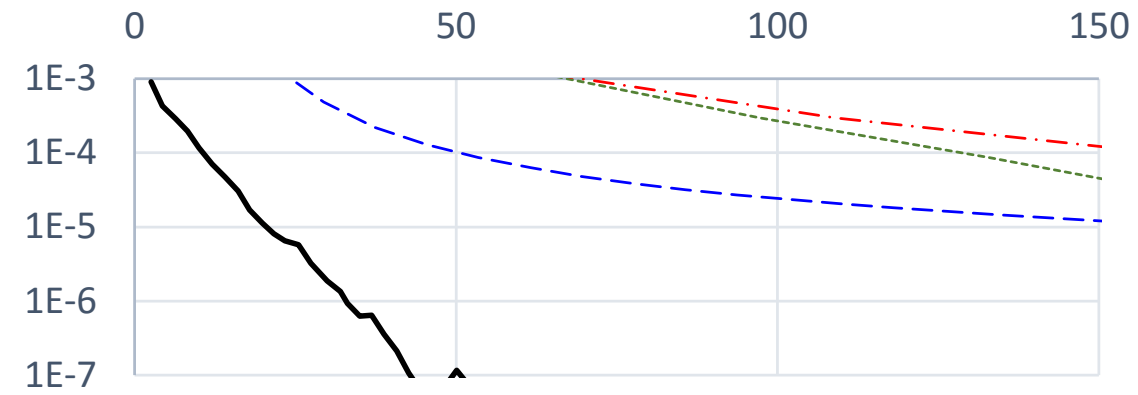

(i) amazon, $k = 30$, `Fnorm` metric

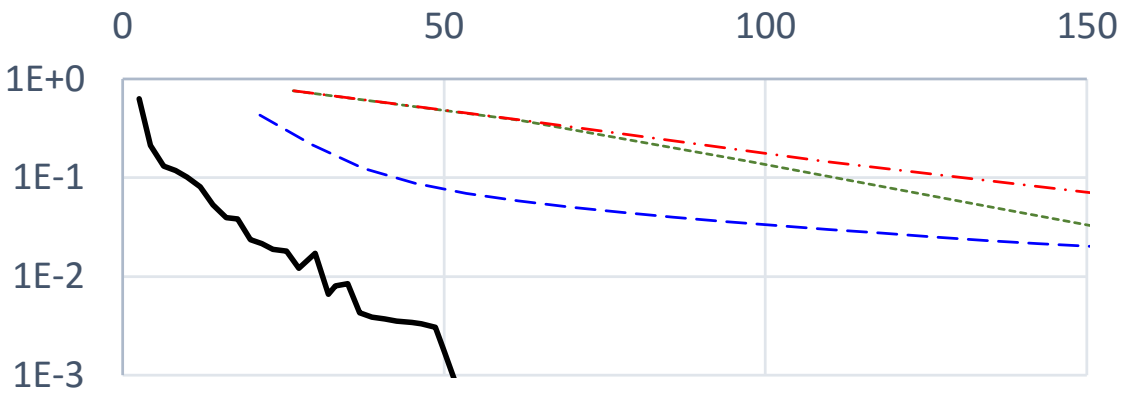

(j) amazon, $k = 30$, `rayleigh` metric

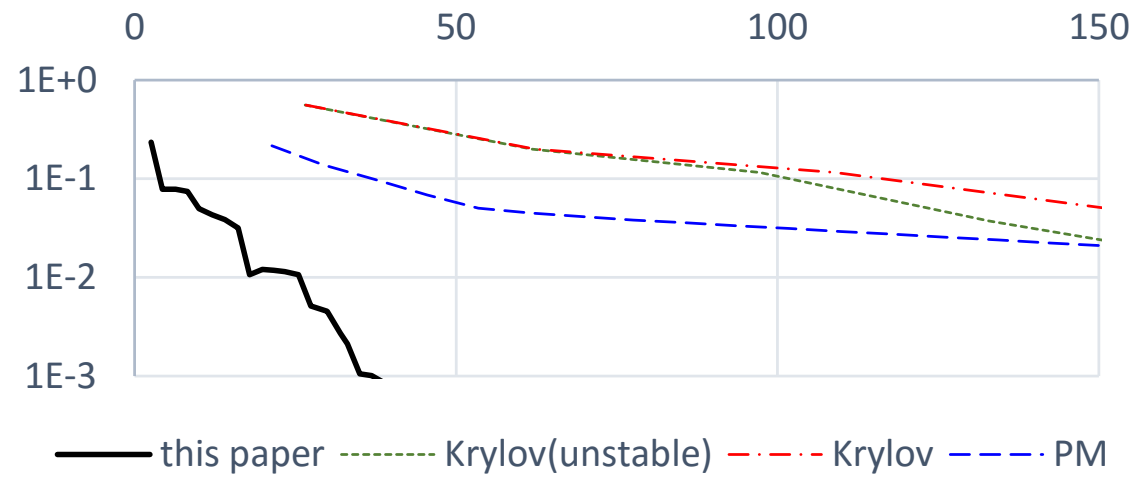

(k) amazon, $k = 30$, `spectral` metric

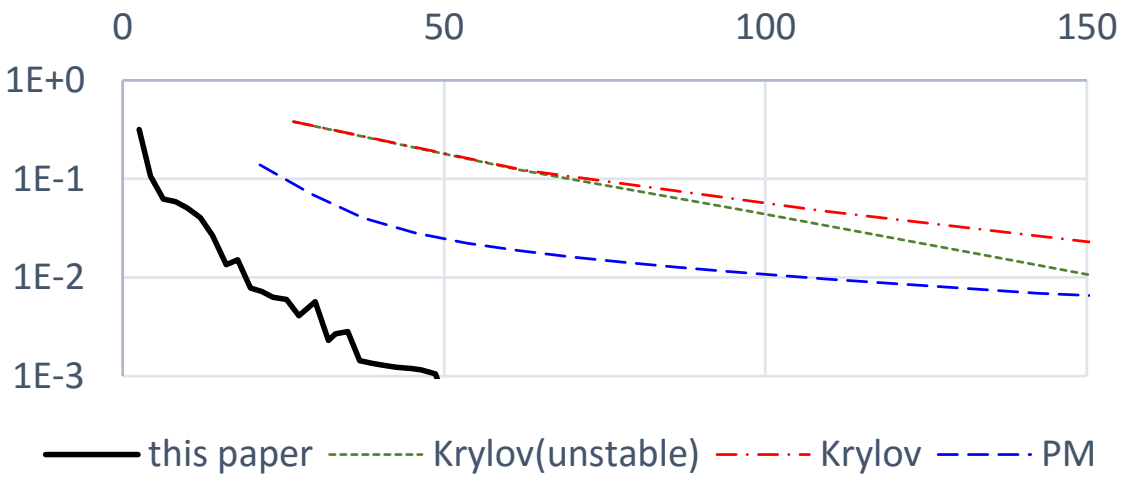

(l) amazon, $k = 30$, `rayleigh(last)` metric

Figure 3: Performance on dataset `amazon`. Relative error ($y$-axis) vs. running time ($x$-axis).

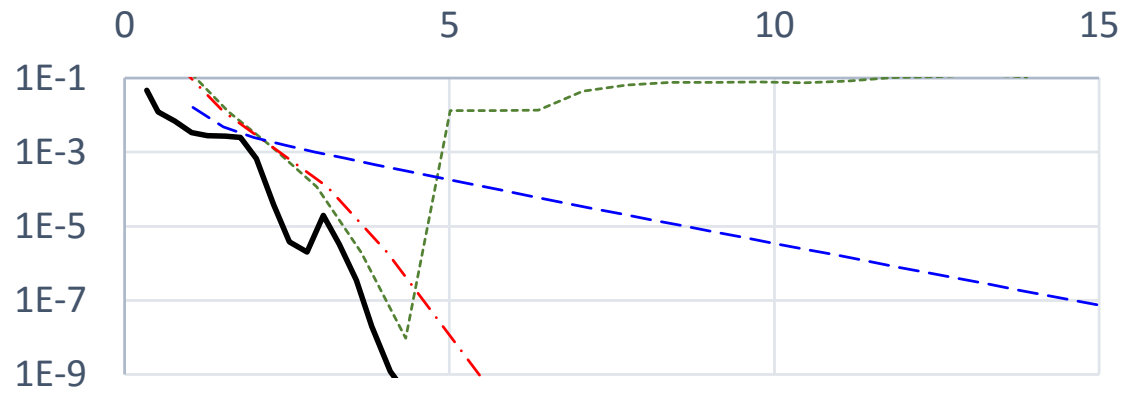

(a) news, $k = 10$, `Fnorm` metric

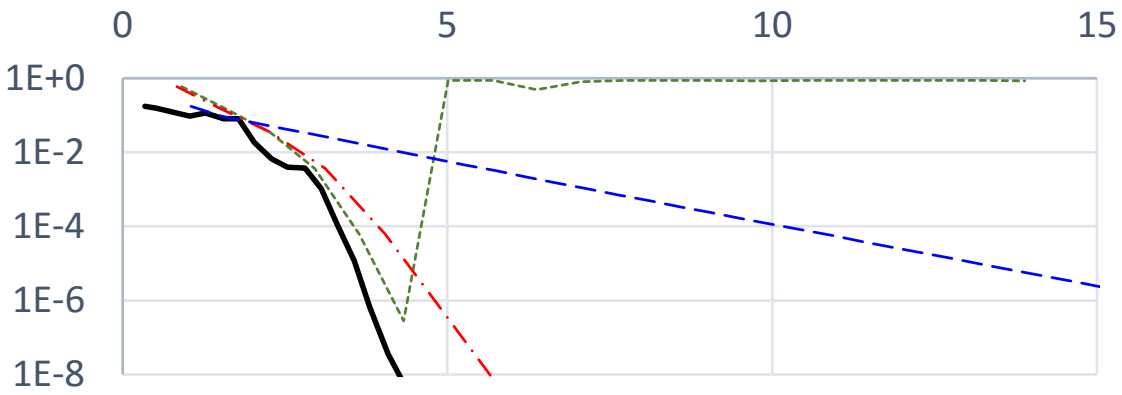

(b) news, $k = 10$, `rayleigh` metric

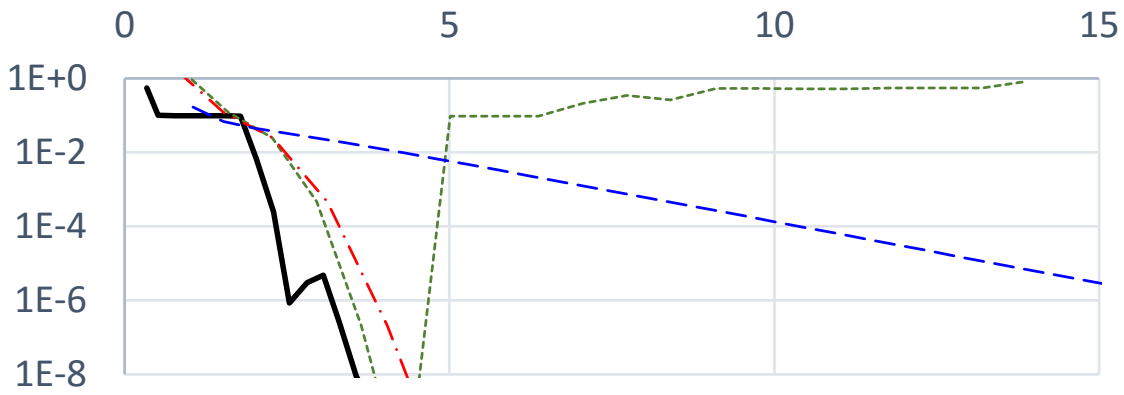

(c) news, $k = 10$, `spectral` metric

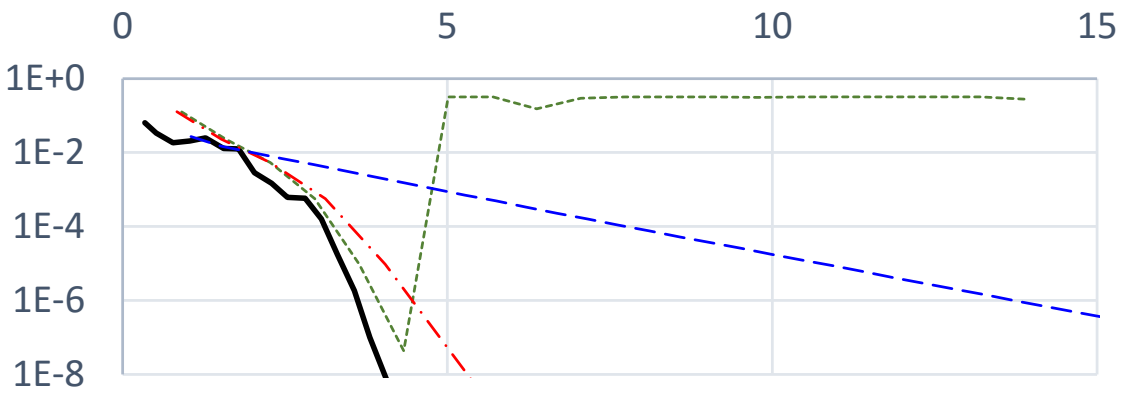

(d) news, $k = 10$, `rayleigh(last)` metric

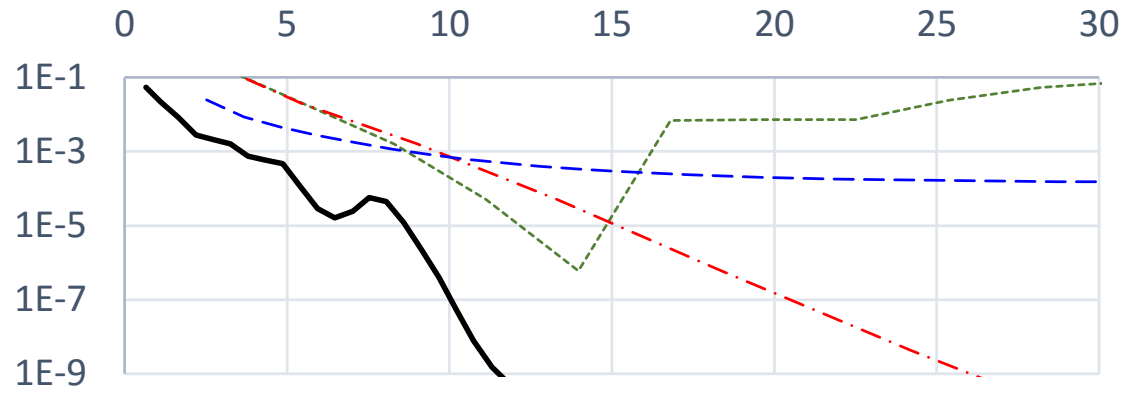

(e) news, $k = 20$, `Fnorm` metric

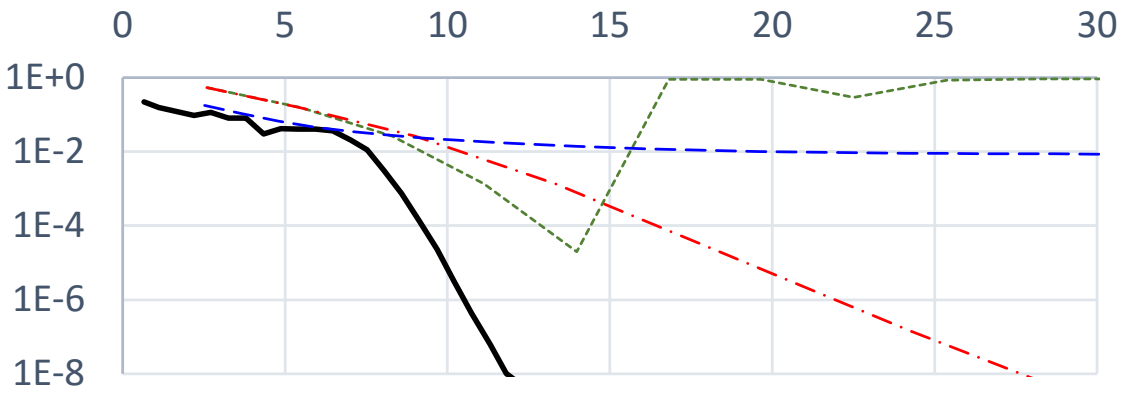

(f) news, $k = 20$, `rayleigh` metric

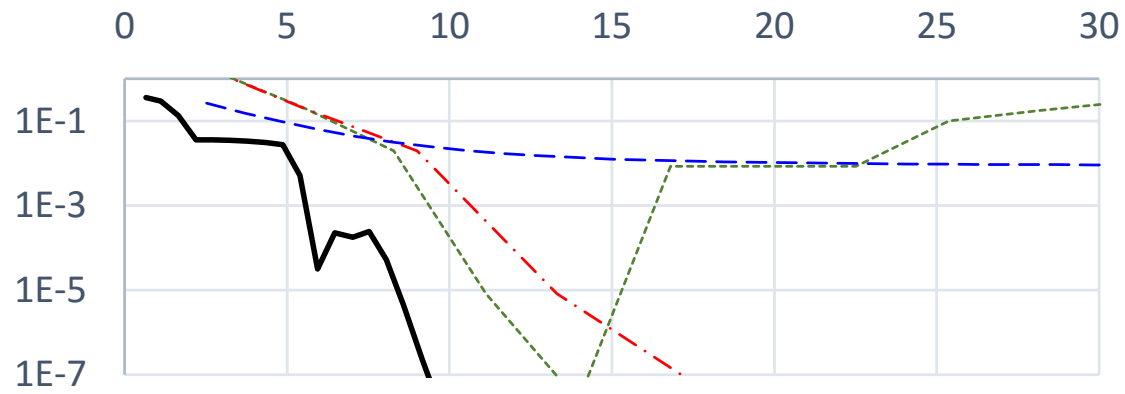

(g) news, $k = 20$, `spectral` metric

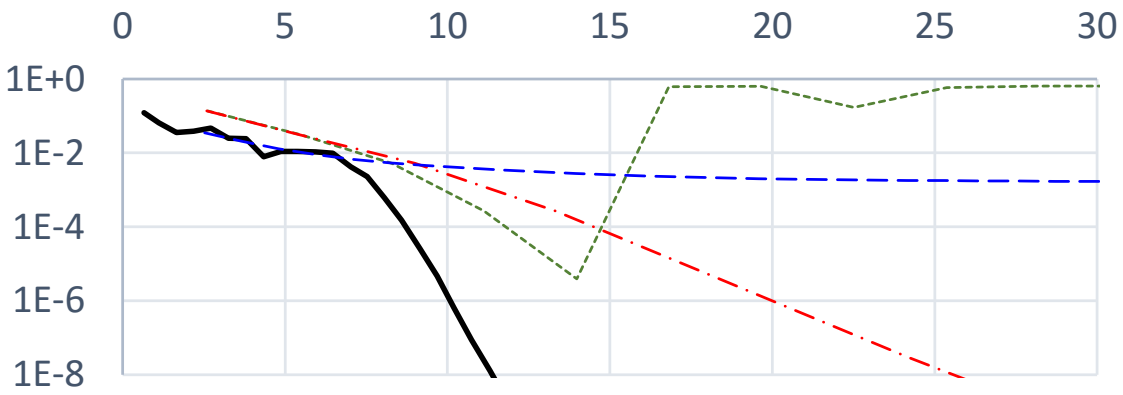

(h) news, $k = 20$, `rayleigh(last)` metric

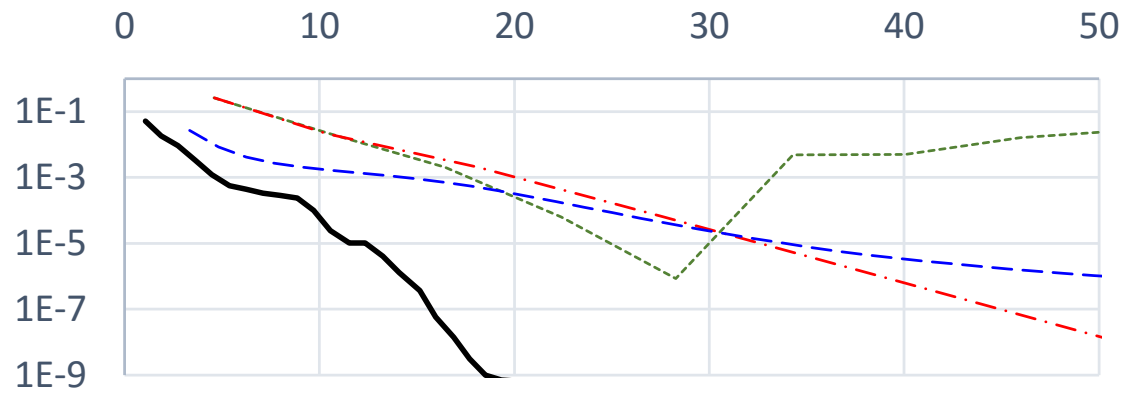

(i) news, $k = 30$, `Fnorm` metric

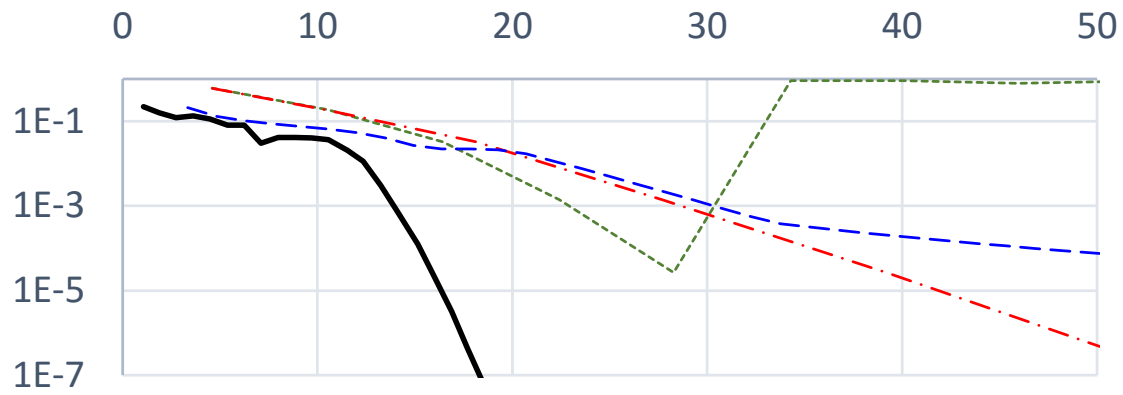

(j) news, $k = 30$, `rayleigh` metric

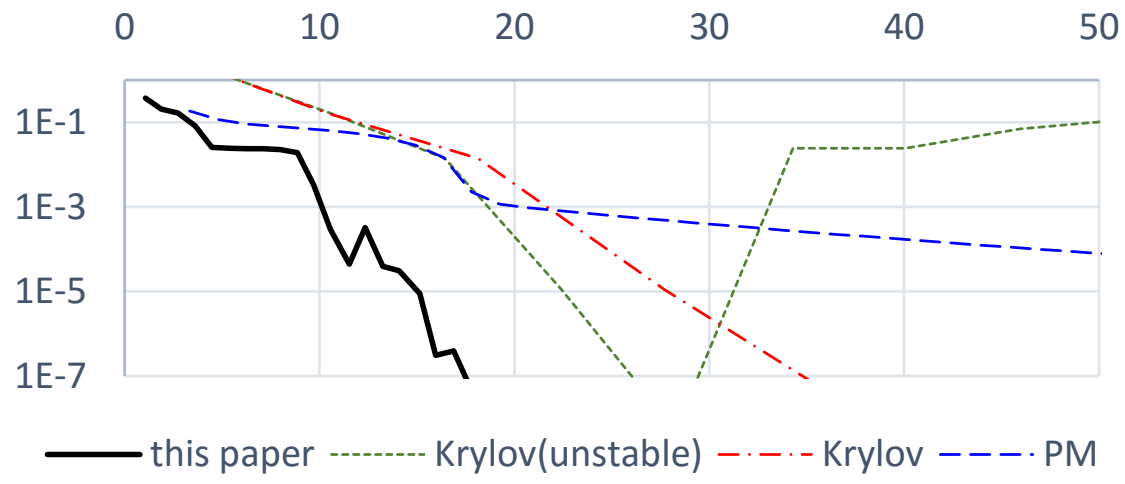

(k) news, $k = 30$, `spectral` metric

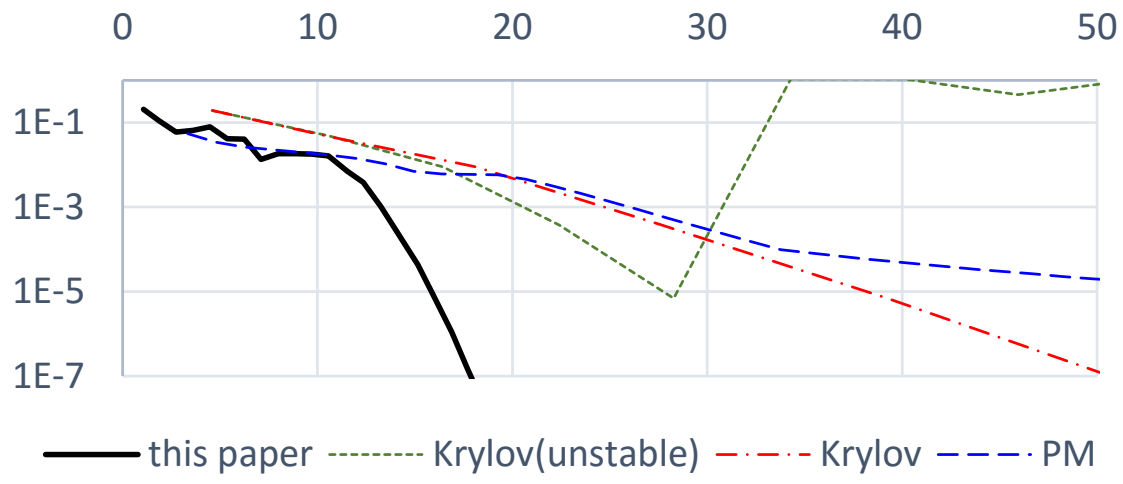

(l) news, $k = 30$, `rayleigh(last)` metric

Figure 4: Performance on dataset `news20`. Relative error ($y$-axis) vs. running time ($x$-axis).

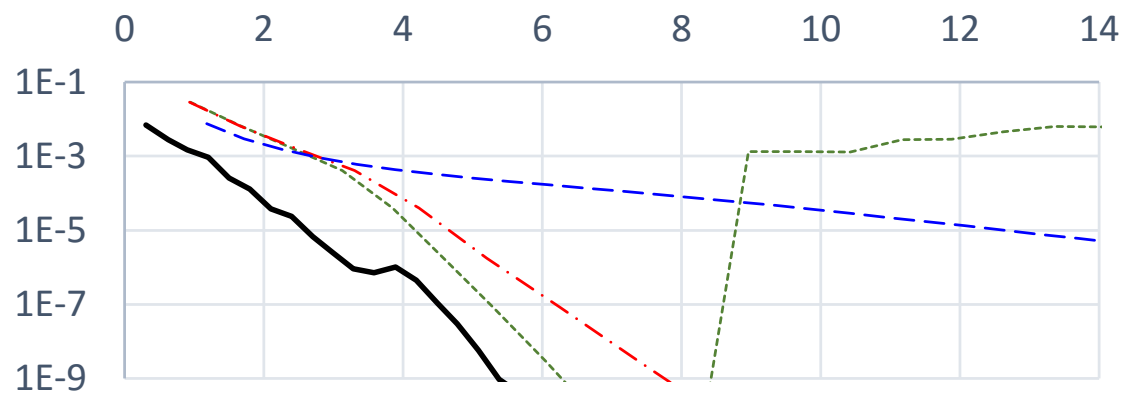

(a) rcv1, $k = 10$, `Fnorm` metric

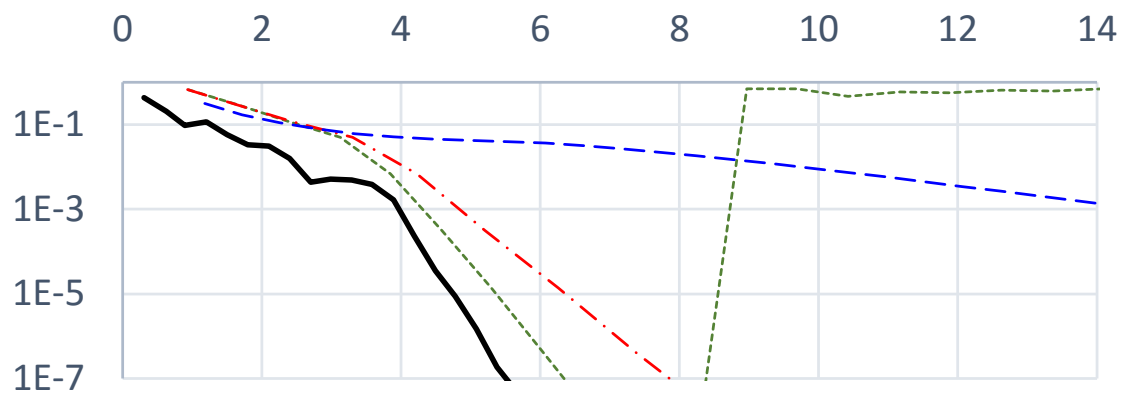

(b) rcv1, $k = 10$, `rayleigh` metric

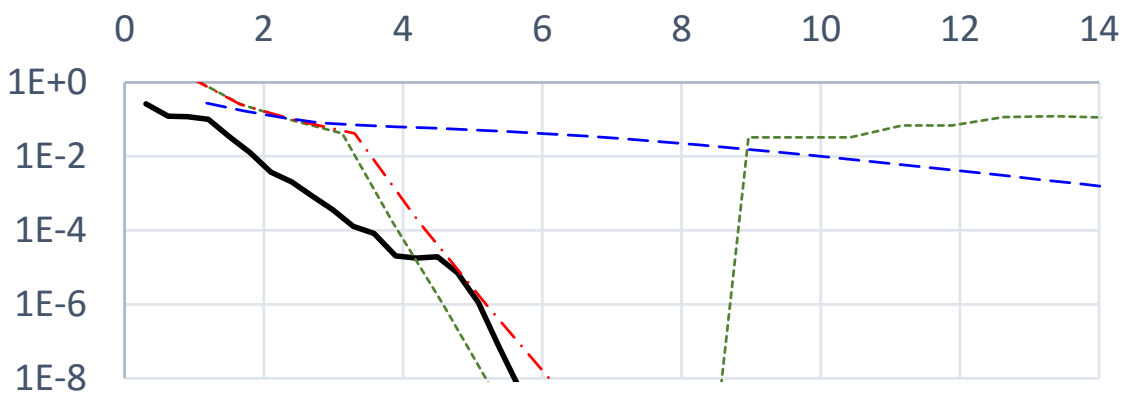

(c) rcv1, $k = 10$, `spectral` metric

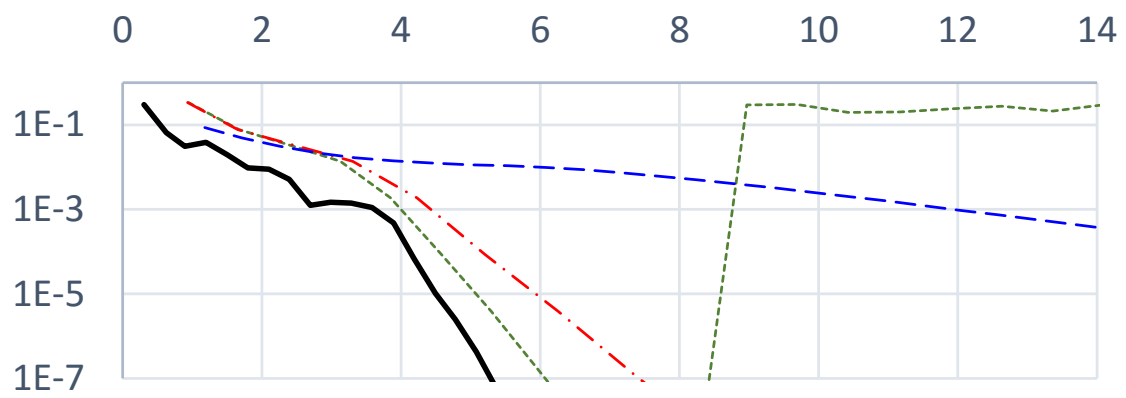

(d) rcv1, $k = 10$, `rayleigh(last)` metric

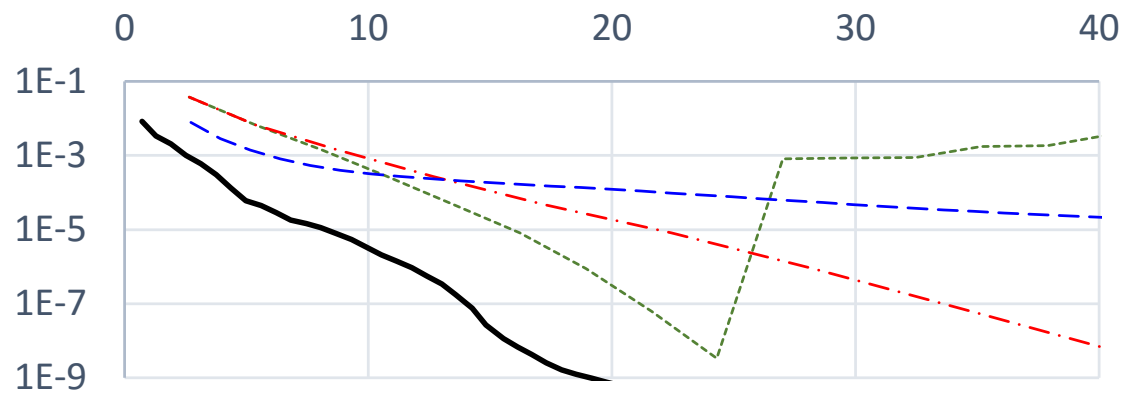

(e) rcv1, $k = 20$, `Fnorm` metric

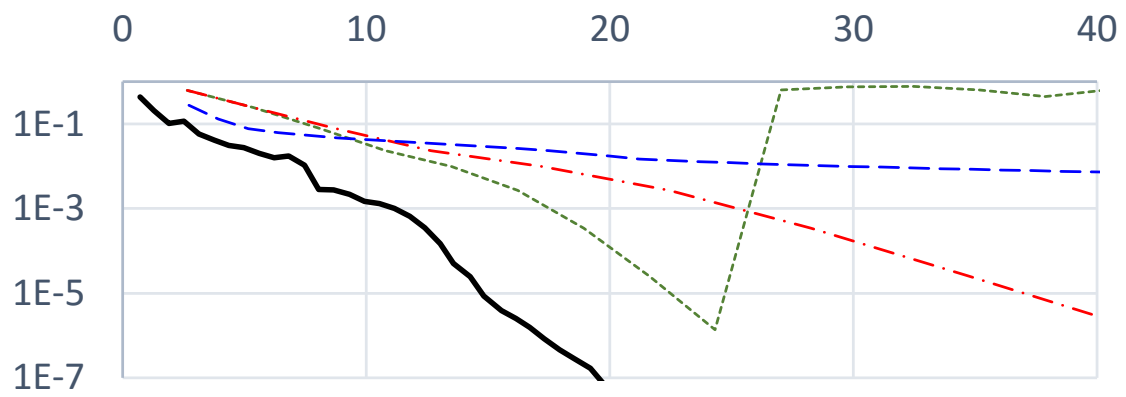

(f) rcv1, $k = 20$, `rayleigh` metric

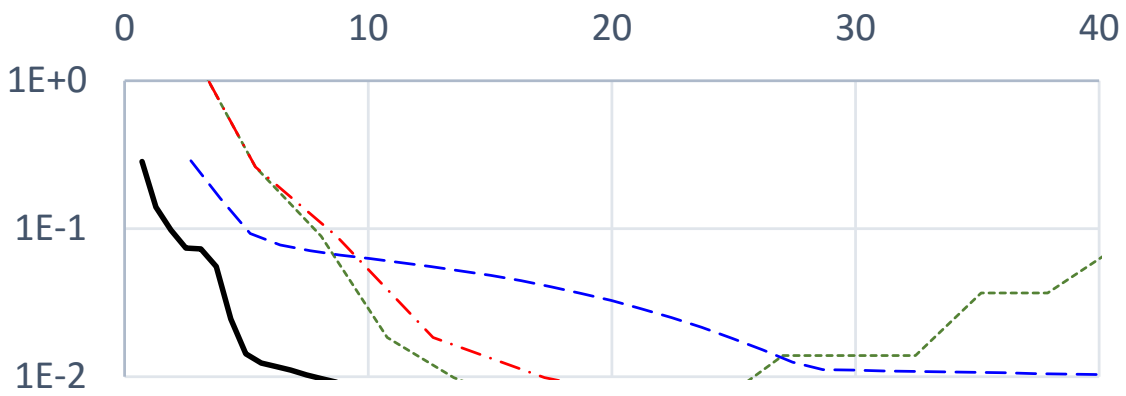

(g) rcv1, $k = 20$, `spectral` metric

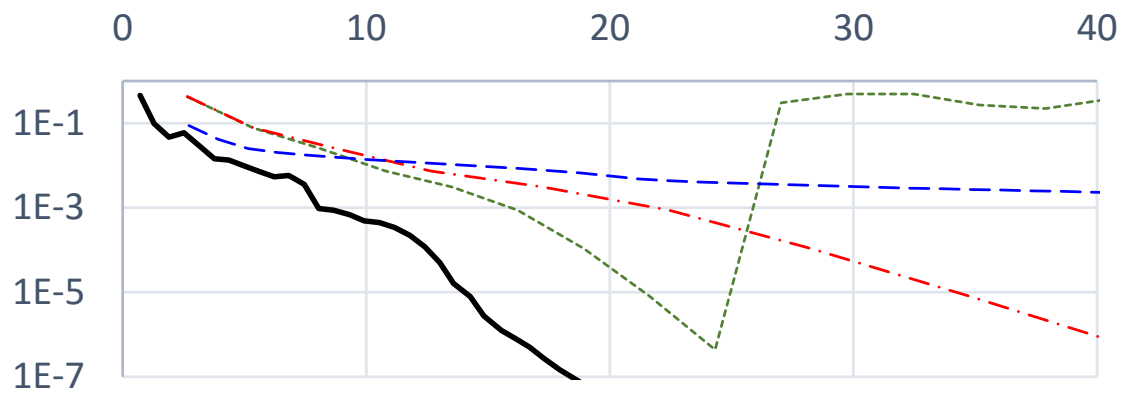

(h) rcv1, $k = 20$, `rayleigh(last)` metric

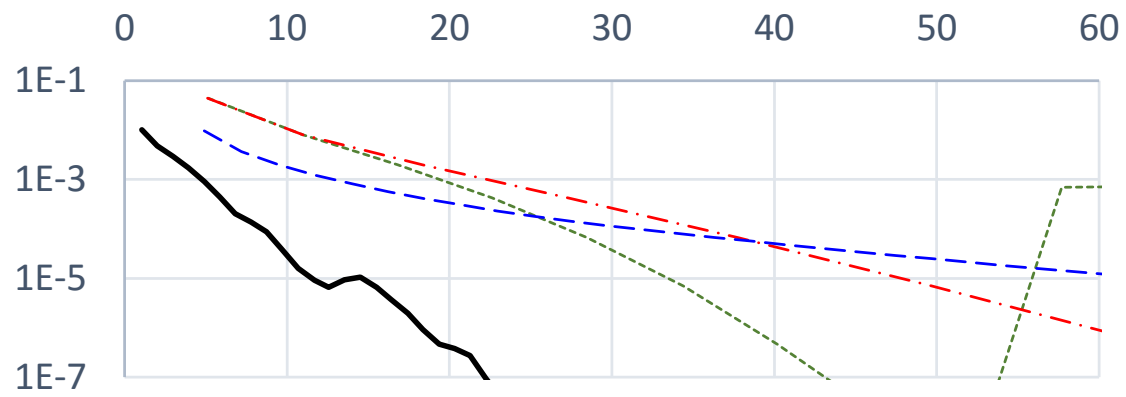

(i) rcv1, $k = 30$, `Fnorm` metric

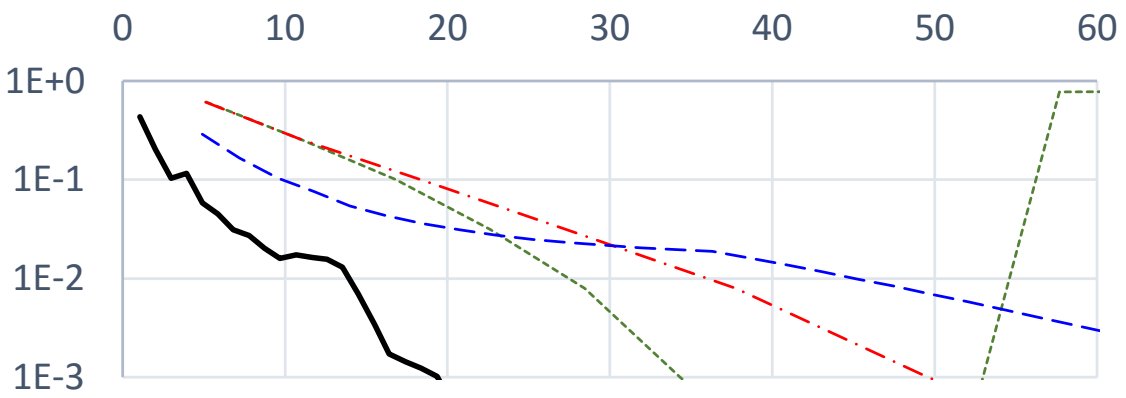

(j) rcv1, $k = 30$, `rayleigh` metric

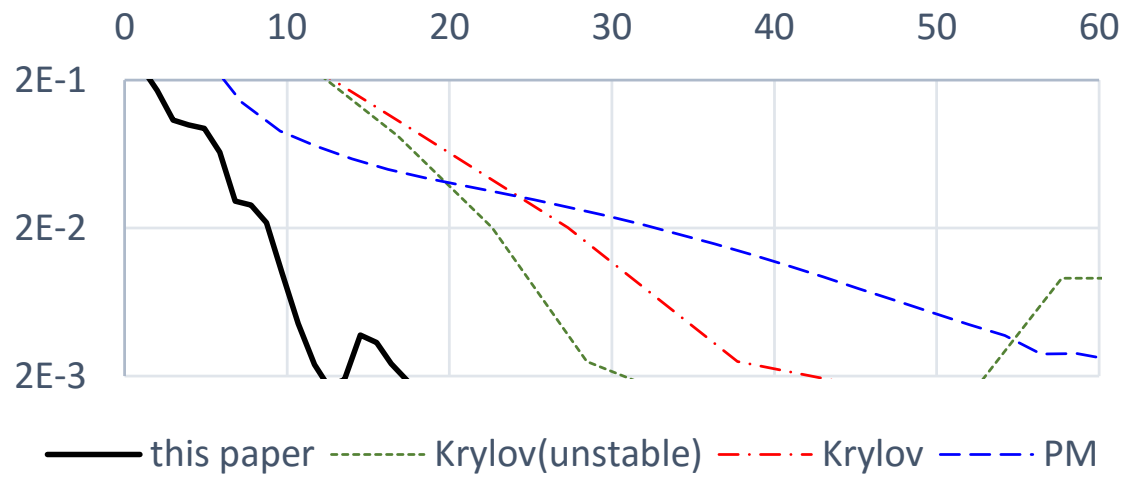

(k) rcv1, $k = 30$, `spectral` metric

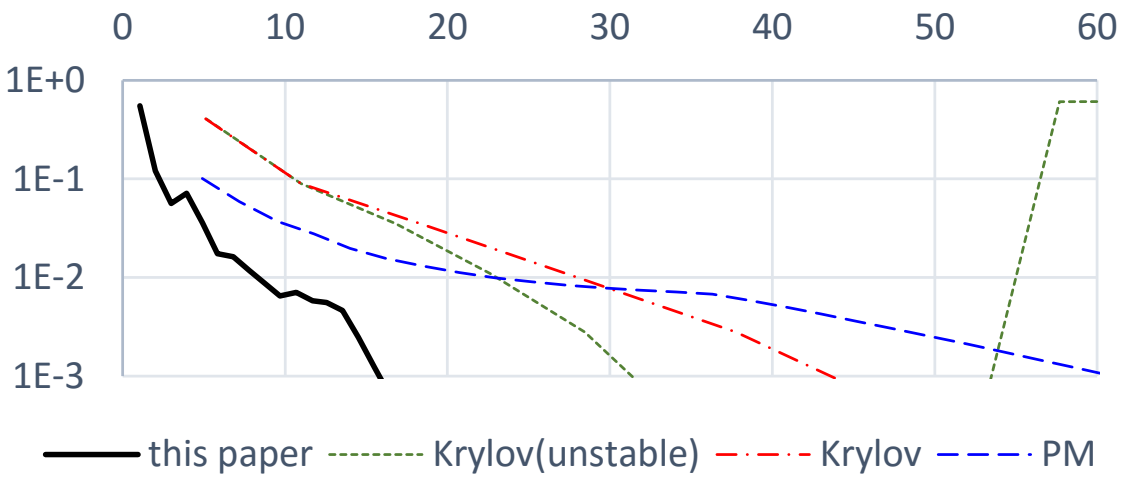

(l) rcv1, $k = 30$, `rayleigh(last)` metric

Figure 5: Performance on dataset `rcv1`. Relative error ($y$-axis) vs. running time ($x$-axis).